\documentclass[11pt,a4paper]{article}
\usepackage[margin=24mm,headheight=14pt]{geometry}
\usepackage[T1]{fontenc}
\usepackage{lmodern,microtype,amsmath,amssymb,amsthm,mathtools,bm,booktabs,array,mathrsfs}
\usepackage[hidelinks]{hyperref}
\usepackage{tikz}
\usepackage{float}
\usepackage{fancyhdr}
\allowdisplaybreaks[1]
\newtheorem{theorem}{Theorem}[section]
\newtheorem{lemma}[theorem]{Lemma}
\newtheorem{proposition}[theorem]{Proposition}
\newtheorem{corollary}[theorem]{Corollary}
\theoremstyle{definition}\newtheorem{definition}[theorem]{Definition}
\theoremstyle{remark}\newtheorem{remark}[theorem]{Remark}
\newtheorem{example}[theorem]{Example}
\newcommand{\R}{\mathbb R}\newcommand{\N}{\mathbb N_0}
\newcommand{\Lo}{\text{\L{}}}\newcommand{\bd}{\mathbf d}
\newcommand{\HH}{\mathrm H}\newcommand{\QQ}{\mathcal Q}\newcommand{\TT}{\mathcal T}
\newcommand{\MM}{\mathcal M}\newcommand{\PP}{\mathcal P}\newcommand{\II}{\mathcal I}
\newcommand{\FF}{\mathcal F}\newcommand{\JJ}{\mathcal J}
\newcommand{\CX}{\mathcal{X}}\newcommand{\CY}{\mathcal{Y}}
\newcommand{\biy}{\bm{\mathit{y}}}\newcommand{\biz}{\bm{\mathit{z}}}
\newcommand{\bx}{\mathbf{x}}\newcommand{\bu}{\mathbf{u}}\newcommand{\bt}{\mathbf{t}}
\newcommand{\bv}{\mathbf{v}}
\newcommand{\bc}{\mathbf{c}}\newcommand{\bg}{\mathbf{g}}
\newcommand{\norm}[1]{\lVert#1\rVert}
\newcommand{\cheb}[1]{\norm{#1}_{1,\mathrm{Cheb}}}
\newcommand{\mon}[1]{\norm{#1}_{1,\mathrm{mon}}}
\DeclareMathOperator{\supp}{supp}\DeclareMathOperator{\dist}{dist}
\DeclareMathOperator{\Lip}{Lip}\DeclareMathOperator{\TV}{d_{TV}}
\DeclareMathOperator{\diag}{diag}
\numberwithin{equation}{section}
\title{Convergence rates of moment--SOS hierarchies\\with correlative sparsity via gluing and projection}
\author{Hoang Anh Tran\thanks{Department of Mathematics, National University of Singapore, Singapore 119076. (tranhoanganh@u.nus.edu).},~
Kim-Chuan Toh\thanks{Department of Mathematics, and Institute of Operations Research and Analytics, National University of Singapore, Singapore 119076. (mattohkc@nus.edu.sg).}}
\date{}
\begin{document}
\maketitle
\begin{abstract}
This paper studies convergence rates of moment--SOS hierarchies with correlative sparsity. Following the method proposed in~\cite{tran}, we bound the Hausdorff distance between the set of sparse truncated moment sequences and sparse truncated pseudo-moment sequences. The central step is a quantitative gluing lemma on a clique tree: local moment errors and separator discrepancies in total variation control the error of a joint approximating measure. Over hypercube, polynomial kernels give the convergence rates $O((\log_2r)^3/r^2)$ for the correlative sparse quadratic module and $O(r^{-2})$ for the correlative sparse preordering, where the certificate degree is $2r$. On compact basic semialgebraic sets, slack variable lifting and global atom projection yield $O((\log_2r)^{3\Lo/2}/r^{\Lo})$ and $O(r^{-\Lo})$ with a {\L}ojasiewicz exponent $0<\Lo \le1$. These bounds hold for fixed-degree Hausdorff distances and polynomial objective gaps. We also extend to result to the correlative sparse quadratic-module rate for scalar objectives with clique-local polynomial matrix constraints.
\end{abstract}


\section{Introduction}
We write $\bx=(x_1,\ldots,x_n)$ and consider the problem of minimizing a polynomial $f\in\R[\bx]$ over a nonempty compact basic semialgebraic set $\CX\subseteq\R^n$:
\begin{equation}\tag{POP}\label{eq:POP}
 f_{\min,\CX}=\min_{\bx\in\CX}f(\bx).
\end{equation}
The feasible set $\CX$ is defined either by polynomial inequalities and equalities as follows:
\begin{equation}\label{eq:scalarset}
 \CX:=\{\bx\in\R^n:g_j(\bx)\ge0\ \forall j\in[m_g],\quad
                       h_\ell(\bx)=0\ \forall\ell\in[m_h]\},
\end{equation}
where $m_g,m_h\ge0$ are the numbers of inequalities and equalities, respectively, and $g_j,h_\ell\in\R[\bx]$, or by a polynomial matrix inequality of the form
\begin{equation}\label{eq:matrixPOP}
 \CX:=\{\bx\in\R^n:G(\bx)\succeq0\},
\end{equation}
where $G\in\mathbb S^p[\bx]$ is a symmetric polynomial matrix of size $p$. Here $[q]=\{1,\ldots,q\}$ and $[0]=\varnothing$, so either constraint family in~\eqref{eq:scalarset} may be absent. We refer to~\eqref{eq:scalarset} and~\eqref{eq:matrixPOP} as the scalar and matrix settings, respectively; the objective $f$ is scalar in both settings. We treat these descriptions separately in the below sections.

For scalar constraints, the moment--sum-of-squares (SOS) hierarchy uses truncated quadratic modules or preorderings as positivity certificates; see Lasserre~\cite{lasserre2001}, Parrilo~\cite{parrilo}, Putinar~\cite{putinar}, and Schm\"udgen~\cite{schmudgen}. Polynomial matrix constraints are treated through polynomial vector multipliers and matrix localizing matrices, as in Scherer and Hol~\cite{schererhol} and Tran and Toh~\cite{tranmatrix}. In both settings considered here, the objective is scalar and the moment formulation optimizes over scalar truncated pseudo-moment sequences subject to the positive semidefiniteness conditions on the moment matrix and the localizing matrices, which results in an SDP (semidefinite programming) relaxation. In the dense case, a dense moment matrix of order $r$ has size $\binom{n+r}{n}$, which limits the orders that can be used in computation. We study the corresponding correlative sparse quadratic module and correlative sparse preordering, and a matrix extension with an explicitly specified localizing matrix, which is used in practice to reduce the size of the SDP relaxation for the sake of calculation. However, sparsity comes with a trade-off in quality. In this paper, we address the question about the efficiency of the correlative sparse moment-SOS hierarchy by exploring the convergence rate of the lower bounds given by the SDP relaxation in terms of the relaxation order $r$. 

Correlative sparsity reduces these matrix sizes by using the variable interactions in the polynomial data. One writes
\begin{displaymath}
 f(\bx)=\sum_{a=1}^{m}f_a(\bx_{I_a}),\qquad I_a\subseteq[n],
\end{displaymath}
and assigns each scalar constraint to a set $I_a$ containing all its variables. In the matrix setting, we use a prescribed family of clique-local polynomial matrices $G_{aj}$; their inequalities can be collected into~\eqref{eq:matrixPOP} by taking $G=\diag(G_{aj})$. Each block $G_{aj}$ is assigned in its entirety to one clique, including the variables in every entry. The dense moment matrix is replaced by matrices on these variable subsets, with equality of the moments on their intersections. The largest moment matrix then has size at most $\binom{w+r}{w}$, where $w=\max_a|I_a|$. A clique tree satisfying the running-intersection property permits compatible local measures to be combined into a measure on all variables.

Following Tran and Toh~\cite{tran,extension}, we quantify the approximation by a Hausdorff distance between the set of truncated moment sequences and the set of truncated pseudo-moment sequences. For a fixed degree $k$, let $A_k$ contain the clique-supported exponents and let $\bv_{A_k}(\bx)=(\bx^\alpha)_{\alpha\in A_k}$. In the notation of Section~\ref{sec:momentsets}, $\MM_k^{\II}(\CX)$ consists of the vectors $\int_\CX\bv_{A_k}\,d\mu$ for probability measures $\mu$ on $\CX$, whereas $\MM_{k,r}^{\II,\mathcal C}(\CX)$ consists of the restrictions to $A_k$ of order-$r$ sparse pseudo-moment sequences. These objects and their clique restrictions are defined in Definitions~\ref{def:sparsepseudo} and~\ref{def:cliquemoments}. For scalar constraints we use $\mathcal C\in\{\QQ,\TT\}$; Section~\ref{sec:matrix} defines the corresponding feasible sets and distances for matrix constraints. Since exact sparse moment sequences extend to feasible sparse pseudo-moment sequences, the Hausdorff distance in~\eqref{eq:distance} reduces to
\begin{equation}\label{eq:introhausdorff}
 \bd_{k,r}^{\mathcal C}(\CX)
 =\sup_{\biy\in\FF_r^{\mathcal C}(\CX)}
   \inf_{\biz\in\MM_k^{\II}(\CX)}\norm{\biy_{A_k}-\biz}_\infty.
\end{equation}
Thus a uniform approximation of the truncated coordinates of every sparse pseudo-moment sequence by exact sparse moments controls the whole relaxation. If $f=\sum_a f_a$ with $\deg f_a\le k$, the coefficient bound recalled in Section~\ref{sec:momentsets} gives
\begin{equation}\label{eq:introgap}
 0\le f_{\min,\CX}-p_r^{\mathcal C}(f;\CX)
 \le F_1\bd_{k,r}^{\mathcal C}(\CX),\qquad
 F_1=\sum_a\mon{f_a}.
\end{equation}
This is the sparse form of the moment-distance estimate of Tran and Toh~\cite[Lemma~2.4]{tran}; see also~\cite{extension}. A bound $\bd_{k,r}^{\mathcal C}(\CX)=O(\psi(r))$ therefore gives the same rate for every fixed-degree clique-supported objective with bounded $F_1$.

The quantitative question is whether this construction preserves the convergence rate available for the dense hierarchy. At a finite order, consistency of clique pseudo-moment sequences defines a sparse functional, whereas a dense positive extension requires additional conditions; see Appendix~\ref{sec:gluingdistinction}. Moreover, independently chosen approximating measures may have different separator marginals, even when the underlying pseudo-moments agree. We address this issue by constructing local measures with controlled marginal discrepancies and proving a quantitative gluing lemma. The resulting joint measure approximates the truncated sparse pseudo-moment sequence. For a constrained set, lifting and projection then turn this ambient approximation into exact moments on the feasible set.

\subsection{Related work}
Waki et al.~\cite{waki} introduced sparse SOS relaxations based on the interaction between variables. Lasserre~\cite{lasserre2006} proved convergence under a running-intersection condition and local compactness assumptions. Grimm et al.~\cite{grimm} gave a direct proof of the corresponding sparse positivity representation. These results explain why variable blocks can be used without losing asymptotic convergence when the overlap structure is appropriate.

Korda et al.~\cite{korda} established polynomial convergence rates for hierarchies with correlative sparsity. Their analysis constructs sparse positivity certificates, with degree estimates depending on the clique sizes and the geometry of the local constraint systems. On the hypercube, their correlative sparse preordering estimate gives order $r^{-2/(w+3)}$ after converting coordinatewise degree bounds to total degree for fixed $w$. The construction below gives order $r^{-2}$ for the explicitly defined correlative sparse preordering. On constrained sets, our estimate uses a global feasibility error bound for the intersection of the local sets. Its exponent should not be identified with the clique-local exponents in~\cite{korda} without an additional argument.

Quantitative estimates for dense hierarchies can be obtained from effective positivity certificates or from approximation of moment sequences. Baldi and Mourrain~\cite{baldi} developed effective Putinar bounds together with moment-approximation estimates. Laurent and Slot~\cite[Theorem~2]{laurent} proved a quadratic convergence rate for the Schm\"udgen hierarchy on the hypercube using polynomial kernels. Gribling et al.~\cite[Theorem~3 and Theorem~7]{kernel} constructed squared polynomial kernels with an error of order $(\log_2 r)^3/r^2$ for the quadratic-module hierarchy on the hypercube. We recall the properties of their squared kernel and the Jackson kernel of Laurent and Slot needed for our gluing lemmas. The additional estimates needed here concern the marginals on clique separators.

Our approach continues the moment-sequence method of Tran and Toh~\cite{tran,extension}. In~\cite{tran}, the convergence analysis estimates the Hausdorff distance between truncated pseudo-moment sets and true moment sets. A finite atomic representation and projection of its atoms transfer approximation on a simple ambient set to approximation on the feasible set. The extension in~\cite{extension} uses polynomial slack variables to lift the pseudo-moment sequence before this projection. We recall these steps in Section~\ref{sec:recall} and develop their correlative sparse counterparts. The work of Tran and Toh on polynomial semidefinite programs~\cite{tranmatrix} also connects convergence estimates with matrix feasibility error bounds.

Term sparsity provides a further reduction of semidefinite blocks. Wang et al.~\cite{tssos,chordal} developed TSSOS and Chordal-TSSOS, and Wang et al.~\cite{cstssos} combined term sparsity with correlative sparsity. Under the appropriate exact-completion conditions, the combined hierarchy recovers the correlative sparse value at a fixed relaxation order. The present paper studies the rate of the correlative sparse hierarchy itself; it does not estimate the additional gap at an incomplete term-sparsity step.

Nie et al.~\cite{niequ} study scalar objectives with clique-local polynomial matrix constraints. Their localizers act on scalar moments, so scalar measure gluing remains applicable. Miller et al.~\cite{miller} also study matrix-valued objectives, for which the corresponding correlative sparse hierarchy may fail to converge. Our matrix extension concerns scalar objectives and full matrix localizers within the prescribed variable cliques.

Blomenhofer et al.~\cite{optimal2026} recently announced improved dense degree estimates. These estimates concern fixed polynomials and do not directly provide compatible clique measures. The rates proved here describe uniform approximation of fixed-degree restrictions of sparse pseudo-moment sequences by the gluing-and-projection construction; we do not assert that they are sharp. A pointwise objective estimate with a constant depending on the objective requires a separate uniformity argument before it gives a Hausdorff estimate for a moment set.

\subsection{Contributions}
We establish quantitative convergence of moment--SOS hierarchies with correlative sparsity. Fix a target degree $k$, the polynomial, and a variable clique tree, and let the relaxation order $r$ increase. The certificate degree is $2r$. The results approximate the fixed-degree restriction of every sparse pseudo-moment sequence by an exact sparse moment sequence on $\CX$. In the notation of Section~\ref{sec:momentsets}, we bound $\bd_{k,r}^{\mathcal C}(\CX)$ uniformly over $\MM_{k,r}^{\II,\mathcal C}(\CX)$. Equations~\eqref{eq:introhausdorff}--\eqref{eq:introgap}, following the moment-approximation approach of~\cite{tran,extension}, convert these Hausdorff bounds into convergence rates for the objective values. The rates are uniform over clique-supported objectives of fixed degree with bounded coefficient norm $F_1$. These results are obtained based on a series of gluing lemmas, where we can glue measures, kernels, and moment sequences with controlled bounds that we present in what follows. 

\paragraph{Quantitative gluing and rate transfer.}
Lemmas~\ref{lem:gluing} and~\ref{lem:compactgluing} establish the quantitative gluing estimates over measures defined on cliques. Consequently, Theorem~\ref{thm:ratetransfer} then separates local dense approximation from separator compatibility. If the local measures have moment error $a_k$ and their separator discrepancies in total variation have maximum tree-path sum $\delta_*$, one joint probability measure has moment error
\begin{displaymath}
 e_k=a_k+\delta_*.
\end{displaymath}
Thus gluing preserves a dense base rate when both errors have that order. Running intersection supplies the conditional-extension construction; Example~\ref{ex:cycle} shows why overlap consistency alone is insufficient on a cyclic clique family. Appendix~\ref{sec:gluingdistinction} recalls the literature distinguishing this measure construction from matrix completion and moment extension.

\paragraph{Rates on the hypercube.}
Corollaries~\ref{thm:boxT} and~\ref{thm:boxQ} give the correlative sparse rates
\begin{displaymath}
 \bd_{k,r}^{\TT}([0,1]^n)=O(r^{-2}),\qquad
 \bd_{k,r}^{\QQ}([0,1]^n)=O((\log_2r)^3/r^2).
\end{displaymath}
The same rates hold for the objective gaps. For the correlative sparse preordering, Lemma~\ref{lem:compatiblekernel} glues the exactly normalized kernel measures without changing their clique marginals. For the correlative sparse quadratic module, Lemma~\ref{lem:squaredseparator} bounds the separator discrepancies in total variation, and Lemma~\ref{lem:gluing} gives a correction of smaller order than the local moment error. Hence neither construction loses a convergence exponent in the gluing step.

\paragraph{Rates on compact sets with scalar constraints.}
The compact gluing lemma also controls constraint violations. Slack variable lifting and global atom projection give
\begin{equation}\label{eq:summarytransfer}
 \bd_{k,r}^{\mathcal C}(\CX)\le e_k+k c_{\Lo}(C_He_t)^{\Lo/2},
 \qquad \mathcal C\in\{\QQ,\TT\},
\end{equation}
where $t$ is a fixed constraint-dependent degree, $C_H$ is the sum of the coefficient norms of the squared lifted equations, and the ambient errors are evaluated at an order proportional to $r$. The scalar \L{}ojasiewicz inequality~\cite[Corollary~2.6.7]{bcr}, used in Tran and Toh~\cite{tran}, gives a constant $c_{\Lo}>0$ and an exponent $0<\Lo\le1$ such that
\begin{equation}\label{eq:LoIntro}
 \dist(\bx,\CX)\le c_{\Lo}\rho(\bx)^{\Lo}\qquad(\bx\in[0,1]^n),
\end{equation}
with $\rho$ defined in~\eqref{eq:residual}. Theorem~\ref{thm:CS} then gives the rates in Table~\ref{tab:rates}. The exponent $\Lo/2$ comes from projection, as in the dense argument. The slack variables preserve the original clique tree and its separators.

\begin{table}[ht]
\centering
\small
\setlength{\tabcolsep}{8pt}
\renewcommand{\arraystretch}{1.25}
\caption{Correlative sparse convergence rates and the Hausdorff distances over the hypercube and compact sets with different positivity certificates.}\label{tab:rates}
\begin{tabular}{@{}lcc@{}}
\toprule
Correlative sparse cone & Hypercube $[0,1]^n$ & Compact set $\CX$\\
\midrule
Quadratic module & $O((\log_2r)^3/r^2)$ & $O((\log_2r)^{3\Lo/2}/r^{\Lo})$\\
 & Corollary~\ref{thm:boxQ} & Theorem~\ref{thm:CS}\\
\addlinespace[3pt]
Preordering & $O(r^{-2})$ & $O(r^{-\Lo})$\\
 & Corollary~\ref{thm:boxT} & Theorem~\ref{thm:CS}\\
\bottomrule
\end{tabular}
\end{table}

\paragraph{Polynomial matrix constraints.}
For scalar objectives with clique-local polynomial matrix inequalities, Theorem~\ref{thm:matrixQ} gives $O((\log_2r)^{3\Lo/2}/r^{\Lo})$ for the correlative sparse matrix quadratic module. Corollary~\ref{cor:matrixB} gives $O(r^{-\Lo})$ for its specified strengthening with box-product localizing matrices. Here $\Lo$ is an exponent furnished by the matrix \L{}ojasiewicz inequality of Dinh and Pham~\cite[Theorem~4.1]{dinhpham}, recalled in Tran and Toh~\cite[Theorem~2.4]{tranmatrix} and applied in~\eqref{eq:matrixerrorbound}. Lemma~\ref{lem:matrixvariance} bounds the squared negative-eigenvalue residual of the local kernel measures. This estimate permits global projection without a polynomial matrix square-root lift.

Compared with the sparse certificate approach of Korda et al.~\cite{korda}, our hypercube proof constructs a joint approximating measure directly, without decomposing a positive objective into positive clique polynomials. It gives the quadratic rate for the correlative sparse preordering with an exponent independent of the clique width. The constants still depend on the clique system. On a constrained set, the exponent $\Lo$ belongs to the global feasible intersection, so we do not claim a uniform improvement over estimates based on different local error bounds. The dense kernel estimates are inputs to the argument; the quantitative separator estimates, compact gluing, and their interaction with lifting and projection are the additional steps.

The compact-set results apply after a coordinatewise affine change of variables and the explicit addition of bounding-hypercube generators in every containing clique. They concern this augmented constraint description. The correlative sparse preordering in Table~\ref{tab:rates} has SOS polynomial multipliers on generator products. It is different from the normalized Krivine--Stengle and extended-Handelman cones with nonnegative scalar coefficients, whose dense base rate in~\cite{extension} is $O(r^{-1})$.

 We present the following outline of the paper: Section~\ref{sec:setup} gives the necessary notation and recalls the moment--SOS hierarchies, correlative sparsity, their convergence assumptions, and the lifting-and-projection method. Section~\ref{sec:method} proves the gluing lemmas and formulates the rate-transfer principle. The two main applications follow: Section~\ref{sec:box} establishes the hypercube rates, and Section~\ref{sec:projection} establishes the rates on compact sets. Section~\ref{sec:matrix} extends the argument to scalar objectives with polynomial matrix constraints. The appendices contain auxiliary evaluation and duality arguments and the distinction between measure gluing and dense pseudo-moment completion.

\section{Preliminaries}\label{sec:setup}
We fix the notation and recall the moment--SOS hierarchies, their correlative sparse counterparts, and the assumptions that ensure convergence. We also recall the moment distances and the lifting-and-projection method used in the rate estimates. The hierarchy convergence results are recalled with references to the literature. We include a proof of the atomic projection estimate, Lemma~\ref{lem:atomic}, because it is used in both the scalar and matrix applications.

\subsection{Notation}
We write $\N=\{0,1,2,\ldots\}$, $[n]=\{1,\ldots,n\}$, and $|\alpha|=\sum_i\alpha_i$ for $\alpha\in\N^n$. The support of an exponent is the set of its nonzero coordinates; the support of a polynomial is the set of exponents with nonzero coefficients. The spaces $\R[\bx]_k$ and $\Sigma[\bx]$ denote polynomials of total degree at most $k$ and sums of squares, respectively. For $p=\sum_\alpha p_\alpha \bx^\alpha$, put $\mon p=\sum_\alpha|p_\alpha|$. We use $\PP(\CY)$ for Borel probability measures on a compact set $\CY$. All distances between points are Euclidean unless stated otherwise. We use plain $x,t$ for univariate variables and bold $\bx=(x_1,\ldots,x_n)$, $\bt=(t_1,\ldots,t_n)$ for vectors; their coordinates $x_i,t_i$ remain plain. Thus $\R[x]$ denotes the univariate polynomial ring and $\R[\bx]$ the multivariate one. Other vectors, such as the slack vector $\bu$ and coefficient vector $\bv$, are also bold. Lifted points are written $\widehat{\bx}=(\bx,\bu)$, retaining $\bx$ for the original variables.
Let $\HH^I=[0,1]^I$. For a set of coordinates $J$, write $\pi_J$ for coordinate projection and $(\pi_J)_\#\mu$ for the corresponding marginal.

\subsection{Moment--SOS hierarchies}
If $\biy=(y_\alpha)_{|\alpha|\leq2r}$, its Riesz functional is $\ell_{\biy}(p)=\sum_\alpha p_\alpha y_\alpha$. The moment matrix and a localizing matrix are
\begin{displaymath}
 M_r(\biy)=(y_{\beta+\gamma})_{\beta,\gamma\in B_r},\qquad
 M_{r-d_g}(g \biy)=\left(\sum_\delta g_\delta y_{\delta+\beta+\gamma}\right)_{\beta,\gamma\in B_{r-d_g}},
\end{displaymath}
where $B_r=\{\alpha\in\N^n:|\alpha|\leq r\}$ and $d_g=\lceil\deg g/2\rceil$. Positive semidefiniteness is equivalent to nonnegativity of $\ell_{\biy}$ on the corresponding squared polynomials and localizing products. A pseudo-moment sequence is not assumed to have a representing measure.

We first recall the hierarchy for the scalar description~\eqref{eq:scalarset}. The matrix localizers for~\eqref{eq:matrixPOP} are recalled separately in Section~\ref{sec:matrixrecall}. The dense truncated quadratic module and equality ideal are
\begin{displaymath}
\begin{aligned}
 \QQ_{2r}(\CX)&=\left\{\sum_{j=0}^{m_g}\sigma_jg_j:
 \sigma_j\in\Sigma[\bx],\ \deg(\sigma_jg_j)\le2r\right\},\qquad g_0=1,\\
 \JJ_{2r}&=\left\{\sum_\ell v_\ell h_\ell:
 v_\ell\in\R[\bx],\ \deg(v_\ell h_\ell)\le2r\right\},
\end{aligned}
\end{displaymath}
where $m_g$ is the number of scalar inequalities and the equality sums range over $\ell\in[m_h]$. The preordering $\TT_{2r}(\CX)$ replaces the generators in the first expression by all their square-free products. Each summand still has degree at most $2r$. These cones depend on the chosen generators, not only on $\CX$. For $\mathcal C\in\{\QQ,\TT\}$, the dense moment and SOS values are
\begin{displaymath}
\begin{aligned}
 p_r^{\mathrm{den},\mathcal C}(f;\CX)
 &=\inf\{\ell_{\biy}(f):y_0=1,\ \ell_{\biy}(\mathcal C_{2r}(\CX))\ge0,
                    \ \ell_{\biy}(\JJ_{2r})=0\},\\
 s_r^{\mathrm{den},\mathcal C}(f;\CX)
 &=\sup\{\lambda:f-\lambda\in\mathcal C_{2r}(\CX)+\JJ_{2r}\}.
\end{aligned}
\end{displaymath}
Evaluation at a probability measure supported on $\CX$ is feasible. Hence both values are lower bounds on $f_{\min,\CX}$. A moment sequence represented by such a measure will be called an exact moment sequence. The sparse relaxations retain only part of the dense positivity conditions or moment coordinates.

The untruncated quadratic module plus ideal is Archimedean if it contains
$R-\sum_i x_i^2$ for some $R>0$. Under this condition and $\CX\ne\varnothing$, the dense moment and SOS bounds increase to $f_{\min,\CX}$, by the moment--SOS construction of Lasserre~\cite{lasserre2001} and Putinar's Positivstellensatz~\cite{putinar}. For the preordering, the corresponding convergence holds on a nonempty compact $\CX$ by Schm\"udgen's Positivstellensatz~\cite{schmudgen}. These hypotheses guarantee asymptotic convergence; finite convergence and equality or attainment of the two values at a given order are separate questions. We next recall the assumptions under which sparsity preserves convergence.

\subsection{Correlative sparsity and convergence assumptions}
We recall the correlative-sparsity construction of Waki et al.~\cite{waki} and the convergence conditions of Lasserre~\cite{lasserre2006}. For a polynomial $p$, write
\begin{displaymath}
 \operatorname{var}(p)=\bigcup_{\alpha\in\supp p}\supp\alpha\subseteq[n].
\end{displaymath}
The correlative sparsity graph has vertex set $[n]$. Distinct vertices $i,j$ are adjacent if they occur together in an objective monomial, or if they both belong to $\operatorname{var}(g_\ell)$ or to $\operatorname{var}(h_\ell)$ for some constraint. In particular, the entire variable set of each constraint is made a clique, even if its variables occur in different monomials.

A graph is chordal if every cycle of length at least four has a chord. A chordal extension is a chordal graph on the same vertices containing all the original edges. Let $\II=(I_a)_{a=1}^m$ be its maximal cliques. The construction ensures that
\begin{equation}\label{eq:csdecomposition}
 f=\sum_{a=1}^m f_a,\qquad f_a\in\R[\bx_{I_a}],
\end{equation}
and that each constraint can be assigned to a clique containing all its variables. The decomposition and assignment need not be unique. We also allow a prescribed family $\II$ covering $[n]$, provided it has these properties and satisfies the running-intersection property.

\begin{definition}[Running-intersection property]\label{def:rip}
The family $\II$ satisfies the running-intersection property if its members can be ordered so that, for every $a>1$, there is $p(a)<a$ with
\begin{equation}\label{eq:rip}
 I_a\cap\bigcup_{b<a}I_b\subseteq I_{p(a)}.
\end{equation}
Equivalently, there is a tree $\mathscr T$ on the clique indices such that the cliques containing any fixed variable form a connected subtree. Such a tree is called a clique tree. For an edge $e=\{a,b\}$, the set $S_e=I_a\cap I_b$ is its separator.
\end{definition}

The maximal cliques of a chordal graph admit a clique tree; disconnected components can be joined by edges with empty separators. Fix a root and write $\omega$ for the largest number of edges on a root-to-node path and $w=\max_a|I_a|$. A chordal extension provides maximal variable cliques with running intersection; these are not two independent assumptions on the chosen clique family. The original data graph may be nonchordal, in which case its extension enlarges the local blocks and the enforced separators. The tree structure is what permits a general marginal-gluing argument~\cite{lasserre2006}. Example~\ref{ex:cycle} shows why pairwise consistency alone is insufficient without it.

\paragraph{Local cones and moment constraints.}
For the scalar description~\eqref{eq:scalarset}, after rescaling a bounding hypercube, write
\begin{displaymath}
 \CX=\{\bx\in[0,1]^n:g_{aj}(\bx_{I_a})\ge0,\ h_{a\ell}(\bx_{I_a})=0
              \text{ for all }a,j,\ell\}\ne\varnothing.
\end{displaymath}
Identically zero constraints are discarded. We explicitly include $b_i=x_i(1-x_i)$ in every clique containing $i$. We take $r\ge1$ with $2r$ at least the degree of the objective and every constraint polynomial. Put
\begin{displaymath}
 B_t^{I_a}=\{\alpha\in\N^{I_a}:|\alpha|\le t\},\qquad
 \bv_t^{I_a}(\bx)=(\bx^\alpha)_{\alpha\in B_t^{I_a}}.
\end{displaymath}
The local quadratic module and ideal at order $r$ are
\begin{displaymath}
\begin{aligned}
 \QQ_{a,2r}(\CX)
 &=\left\{\sigma_0+\sum_{i\in I_a}\sigma_i b_i+\sum_j\sigma'_j g_{aj}:
          \sigma_0,\sigma_i,\sigma'_j\in\Sigma[\bx_{I_a}],
          \ \deg(\text{each summand})\le2r\right\},\\
 \JJ_{a,2r}
 &=\left\{\sum_\ell v_\ell h_{a\ell}:
          v_\ell\in\R[\bx_{I_a}],\ \deg(v_\ell h_{a\ell})\le2r\right\}.
\end{aligned}
\end{displaymath}
The local preordering $\TT_{a,2r}(\CX)$ uses all square-free products of the local inequality generators, with SOS multipliers and the same total-degree bound. Within each summand, products use only the generators assigned to that clique. For $\mathcal C\in\{\QQ,\TT\}$, define
\begin{displaymath}
 \mathcal C_{2r}^{\rm sp}(\CX)=\sum_a\mathcal C_{a,2r}(\CX),\qquad
 \JJ_{2r}^{\rm sp}=\sum_a\JJ_{a,2r}.
\end{displaymath}
We call $\QQ_{2r}^{\rm sp}(\CX)$ the \emph{correlative sparse quadratic module} and $\TT_{2r}^{\rm sp}(\CX)$ the \emph{correlative sparse preordering}, both truncated at degree $2r$. 

For an integer $t\ge0$, let
\begin{equation}\label{eq:sparseexponents}
 A_t=\bigcup_a\{\alpha\in\N^n:\supp\alpha\subseteq I_a,\ |\alpha|\le t\}.
\end{equation}
We identify exponents in $\N^{I_a}$ with their zero extensions to $\N^n$, and identify repeated monomials in this union. Thus $A_{2r}$ contains exactly the moment coordinates used by the order-$r$ correlative sparse relaxation.

For each clique, a truncated real array $\biy^{(a)}=(y_\alpha^{(a)})_{\alpha\in B_{2r}^{I_a}}$ defines a Riesz functional on $\R[\bx_{I_a}]_{2r}$. The local feasibility and separator consistency conditions are
\begin{align}
 \ell_{\biy^{(a)}}(1)&=1,&
 \ell_{\biy^{(a)}}(\mathcal C_{a,2r}(\CX))&\ge0,&
 \ell_{\biy^{(a)}}(\JJ_{a,2r})&=0,\label{eq:feas}\\
 \ell_{\biy^{(a)}}(p)&=\ell_{\biy^{(b)}}(p)
 &&\bigl(p\in\R[\bx_{S_e}]_{2r},\ e=\{a,b\}\bigr).\label{eq:consistency}
\end{align}
For the correlative sparse quadratic module, the inequalities in~\eqref{eq:feas} are precisely
\begin{displaymath}
 M_r(\biy^{(a)})\succeq0,\qquad
 M_{r-1}(b_i\biy^{(a)})\succeq0,\qquad
 M_{r-d_{aj}}(g_{aj}\biy^{(a)})\succeq0,
 \quad d_{aj}=\lceil\deg g_{aj}/2\rceil.
\end{displaymath}
The last matrices are indexed by $B_{r-d_{aj}}^{I_a}$ and have entries
$\sum_\delta(g_{aj})_\delta y_{\delta+\beta+\gamma}^{(a)}$.
Equality constraints mean
$\ell_{\biy^{(a)}}(\bx^\beta h_{a\ell})=0$ whenever
$|\beta|+\deg h_{a\ell}\le2r$. For the correlative sparse preordering, one additionally imposes the localizing matrix for each admissible generator product.

\begin{definition}[Sparse pseudo-moment sequences]\label{def:sparsepseudo}
An order-$r$ \emph{sparse pseudo-moment sequence} for the chosen clique family, constraint description, and cone $\mathcal C\in\{\QQ,\TT\}$ is a truncated array
\begin{displaymath}
 \biy=(y_\alpha)_{\alpha\in A_{2r}}
\end{displaymath}
whose clique restrictions $\biy^{(a)}=(y_\alpha)_{\alpha\in B_{2r}^{I_a}}$ satisfy~\eqref{eq:feas}. We denote the set of these sequences by $\FF_r^{\mathcal C}(\CX)$ and call each restriction $\biy^{(a)}$ a \emph{clique pseudo-moment sequence}. Equivalently, $\biy$ is specified by a family $(\biy^{(a)})_{a=1}^m$ satisfying both~\eqref{eq:feas} and~\eqref{eq:consistency}. The order is $r$, whereas the truncation degree is $2r$.
\end{definition}

Running intersection makes edge consistency in~\eqref{eq:consistency} identify every shared monomial, so the family formulation and the single array indexed by $A_{2r}$ are equivalent; see Lasserre~\cite{lasserre2006}. This equivalence concerns the sparse coordinates only. It does not provide a dense pseudo-moment extension or a representing measure; see Appendix~\ref{sec:gluingdistinction}. Agreement only through a smaller target degree is not the hierarchy used here.

The associated sparse Riesz functional is defined on $E_{2r}=\sum_a\R[\bx_{I_a}]_{2r}$ by $\ell_{\biy}(\bx^\alpha)=y_\alpha$ for $\alpha\in A_{2r}$. Its restriction to each clique polynomial space is $\ell_{\biy^{(a)}}$. The moment and SOS bounds are
\begin{displaymath}
\begin{aligned}
 p_r^{\mathcal C}(f;\CX) &=\inf_{\biy\in\FF_r^{\mathcal C}(\CX)}
              \sum_a\ell_{\biy^{(a)}}(f_a),\\
 s_r^{\mathcal C}(f;\CX) &=\sup\{\lambda:f-\lambda\in\mathcal C_{2r}^{\rm sp}(\CX)+\JJ_{2r}^{\rm sp}\}.
\end{aligned}
\end{displaymath}
The objective is independent of the chosen decomposition~\eqref{eq:csdecomposition}. In general, the correlative sparse bound at a fixed order can be weaker than the dense bound.

\paragraph{Conditions for asymptotic convergence.}
Write $\QQ_a=\bigcup_r\QQ_{a,2r}(\CX)$ and $\JJ_a=\bigcup_r\JJ_{a,2r}$. The local Archimedean condition is
\begin{equation}\label{eq:localarch}
 R_a-\sum_{i\in I_a}x_i^2\in\QQ_a+\JJ_a
 \quad\text{for some }R_a>0,\quad a=1,\ldots,m.
\end{equation}
One can enforce it by adding a redundant ball inequality in every clique, as in Lasserre~\cite[Assumption~3.1]{lasserre2006}. The explicit box generators used in this paper also enforce it. Compactness of the global feasible set, or an Archimedean certificate using all variables together, does not alone establish~\eqref{eq:localarch} for the chosen local modules. The running intersection property and Archimedean condition are sufficient to guarantee the asymptotic convergence presented in the next theorem. 

\begin{theorem}\label{thm:csrecall}
Suppose $\CX\ne\varnothing$, the objective and constraints admit the clique assignments above, $\II$ satisfies~\eqref{eq:rip}, and each local module satisfies~\eqref{eq:localarch}. Then
\begin{equation}\label{eq:csconvergence}
 s_r^{\QQ}(f;\CX)\uparrow f_{\min,\CX},
 \qquad p_r^{\QQ}(f;\CX)\uparrow f_{\min,\CX}.
\end{equation}
More precisely, every strictly positive polynomial
$p\in\sum_a\R[\bx_{I_a}]$ on $\CX$ belongs to
$\sum_a(\QQ_a+\JJ_a)$. The same convergence holds for the correlative sparse preordering hierarchy under these assumptions.
\end{theorem}

For the moment convergence see Lasserre~\cite[Theorem~3.6]{lasserre2006}; for sparse positivity without a nonempty-interior assumption see Grimm et al.~\cite[Theorem~5]{grimm}. Equalities may be included by replacing $h=0$ with $h\ge0$ and $-h\ge0$; we use the full truncated ideal convention in~\eqref{eq:feas}. Furthermore, strong duality also holds due to the local Archimedean condition, which is recorded in Proposition~\ref{prop:duality}.

\subsection{Truncated moment sets and Hausdorff distances}\label{sec:momentsets}
For $k\ge1$, use the exponent set $A_k$ from~\eqref{eq:sparseexponents} and write $\bv_{A_k}(\bx)=(\bx^\alpha)_{\alpha\in A_k}$ for the truncated monomial vector. The following definition defines clique moment sequences and exact sparse moment sequences. 

\begin{definition}\label{def:cliquemoments}
Let $K_a\subseteq[0,1]^{I_a}$ be compact. A \emph{truncated clique moment sequence} of degree $k$ on $K_a$ is an array $\biz^{(a)}=(z_\alpha^{(a)})_{\alpha\in B_k^{I_a}}$ having a representing probability measure $\mu_a\in\PP(K_a)$:
\begin{displaymath}
 z_\alpha^{(a)}=\int_{K_a}\bx_{I_a}^\alpha\,d\mu_a(\bx_{I_a})
 \quad(\alpha\in B_k^{I_a}).
\end{displaymath}
A \emph{truncated exact sparse moment sequence} of degree $k$ on $\CX$ is an array $\biz=(z_\alpha)_{\alpha\in A_k}$ represented by one probability measure $\mu\in\PP(\CX)$, that is,
\begin{displaymath}
 \biz=\int_\CX\bv_{A_k}(\bx)\,d\mu(\bx).
\end{displaymath}
Its clique restrictions are the clique moment sequences of the marginals $(\pi_{I_a})_\#\mu$ on $\pi_{I_a}(\CX)$.
\end{definition}

We omit the word \emph{truncated} when the degree is specified. A clique pseudo-moment sequence is called a clique moment sequence only when a representing measure exists. Agreement of truncated separator moments alone does not identify the separator marginal measures; see Appendix~\ref{sec:gluingdistinction} for the relevant moment-extension and measure-gluing results. Equality of the separator \emph{marginals as measures}, together with running intersection, is sufficient for exact gluing, as recalled in~\cite[Section~6.1 and Lemma~6.4]{lasserre2006} and proved in Corollary~\ref{cor:exactcompact}.

For nonempty compact subsets $\mathcal A,\mathcal B$ of the same finite-dimensional coordinate space, their Hausdorff distance in the maximum norm is
\begin{displaymath}
 \bd_\infty(\mathcal A,\mathcal B)
 =\max\left\{\sup_{\biy\in\mathcal A}\inf_{\biz\in\mathcal B}\norm{\biy-\biz}_\infty,
              \sup_{\biz\in\mathcal B}\inf_{\biy\in\mathcal A}\norm{\biy-\biz}_\infty\right\}.
\end{displaymath}
We use this viewpoint as in Tran and Toh~\cite{tran,extension}.

When $k\le2r$, let $\biy_{A_k}=(y_\alpha)_{\alpha\in A_k}$ denote the degree-$k$ restriction of a sparse pseudo-moment sequence. Define the exact sparse moment set and the truncated sparse pseudo-moment set, respectively, by
\begin{equation}\label{eq:momentsets}
\begin{aligned}
 \MM_k^{\II}(\CX)&=\left\{\int_\CX\bv_{A_k}(\bx)\,d\mu(\bx):\mu\in\PP(\CX)\right\},\\
 \MM_{k,r}^{\II,\mathcal C}(\CX)&=\{\biy_{A_k}:\biy\in\FF_r^{\mathcal C}(\CX)\}.
\end{aligned}
\end{equation}
Thus, $\MM_{k,r}^{\II,\mathcal C}(\CX)$ plays a role as an outer spectrahedral approximation of $\MM_k^{\II}(\CX)$, and our approximation error is evaluated by
\begin{equation}\label{eq:distance}
\begin{aligned}
 \bd_{k,r}^{\mathcal C}(\CX)
 &=\bd_\infty\bigl(\MM_{k,r}^{\II,\mathcal C}(\CX),\MM_k^{\II}(\CX)\bigr)\\
 &=\sup_{\biy\in\FF_r^{\mathcal C}(\CX)}
       \inf_{\biz\in\MM_k^{\II}(\CX)}\norm{\biy_{A_k}-\biz}_\infty.
\end{aligned}
\end{equation}
 The exact sparse moment set is compact due to compactness of $\CX$, and the explicitly included box generators also make the truncated sparse pseudo-moment set compact; see Proposition~\ref{prop:duality}. For the decomposition~\eqref{eq:csdecomposition}, let $F_1=\sum_a\mon{f_a}$. The coefficient-duality estimate of Tran and Toh~\cite[Lemma~2.4]{tran}, applied to the truncated coordinates, gives
\begin{equation}\label{eq:hausdorffobjective}
 0\le f_{\min,\CX}-p_r^{\mathcal C}(f;\CX)
 \le F_1\bd_{k,r}^{\mathcal C}(\CX)
 \qquad\bigl(\deg f_a\le k\le2r\bigr).
\end{equation}
This recalled estimate transfers any uniform Hausdorff rate to the objective gap. Constants in asymptotic bounds may depend on the clique tree, degrees, coefficient norms, and global error-bound data.

We shall also use shifted Chebyshev polynomials
$\tau_j(x)=T_j(2x-1)$ and
$\tau_\alpha(\bx_I)=\prod_{i\in I}\tau_{\alpha_i}(x_i)$.
For the expansion in this unnormalized basis, define
\begin{displaymath}
 p(\bx_I)=\sum_\alpha p_\alpha^{\mathrm{Cheb}}\tau_\alpha(\bx_I),
 \qquad
 \cheb p:=\sum_\alpha|p_\alpha^{\mathrm{Cheb}}|.
\end{displaymath}
This is the $\ell^1$ norm of the Chebyshev coefficient vector, as in~\cite[Section~2.2]{kernel}, after shifting to $[0,1]^I$. It is not an $L^1$ or $L^2$ norm of the function. For clarity, let
\begin{displaymath}
\begin{gathered}
 d\varpi(x)=\frac{dx}{\pi\sqrt{x(1-x)}},\qquad
 \varpi_I=\bigotimes_{i\in I}\varpi_i,\\[3pt]
 \norm{p}_{L^v(\varpi_I)}=
 \left(\int_{[0,1]^I}|p(\bx_I)|^v\,d\varpi_I(\bx_I)\right)^{1/v}
 \quad(v=1,2).
\end{gathered}
\end{displaymath}
These function norms are distinct from $\cheb p$; in particular, orthogonality gives
$\norm{p}_{L^2(\varpi_I)}^2=\sum_\alpha 2^{-|\supp\alpha|}|p_\alpha^{\mathrm{Cheb}}|^2$.
All Chebyshev norms in the kernel estimates below mean the coefficient $\ell^1$ norm $\cheb{\cdot}$. The evaluation bounds and the finite-order duality facts used later are collected in Appendix~\ref{app:auxiliary}.

For probability measures on a hypercube, define the degree-$k$ moment discrepancy and total-variation distance by
\begin{displaymath}
 D_k(\mu,\nu)=\max_{|\alpha|\le k}
       \left|\int \bx^\alpha\,d(\mu-\nu)\right|,
 \qquad \TV(\mu,\nu)=\sup_B|\mu(B)-\nu(B)|,
\end{displaymath}
where the supremum is over Borel sets. Since $0\le \bx^\alpha\le1$ on the hypercube,
\begin{equation}\label{eq:momentmetrics}
 D_k(\mu,\nu)\le\TV(\mu,\nu).
\end{equation}
We will use total variation only to control separator and clique marginal errors of probability measures.
\subsection{The lifting-and-projection method for moment sequences}\label{sec:recall}
We recall the lifting-and-projection argument of Tran and Toh~\cite{tran,extension}. This method explains how moment approximation on an ambient set yields approximation on the feasible set. Let $\CY$ be a compact ambient set containing $\CX$, and let $A$ be a finite set of exponents. Write $\bv_A(\bx)=(\bx^\alpha)_{\alpha\in A}$. The first step is to approximate a pseudo-moment vector $\biy_A$ by
\begin{displaymath}
\overline{\biy}_A=\int_{\CY}\bv_A(\bx)\,d\mu(\bx),\qquad \mu\in\PP(\CY).
\end{displaymath}
Suppose that the approximation controls both these moments and a nonnegative violation function. A finite atomic representation, which exists according to Tchakaloff's theorem (see Bayer and Teichmann~\cite[Corollary~2]{bayer}), then permits the atoms to be projected individually onto $\CX$ to form an exact truncated moment sequence. The details are presented as follows. 

\paragraph{The semialgebraic \L{}ojasiewicz inequality.}
We recall the inequality used in the moment-approximation argument of Tran and Toh~\cite{tran}. Let $\CY$ be compact and semialgebraic, and let $\varphi,\psi:\CY\to\R$ be continuous semialgebraic functions with $\varphi^{-1}(0)\subseteq\psi^{-1}(0)$. By Bochnak et al.~\cite[Corollary~2.6.7]{bcr}, there exist $c_{\Lo}>0$ and $\Lo>0$ such that
\begin{displaymath}
 |\psi(\bx)|\le c_{\Lo}|\varphi(\bx)|^{\Lo}
 \qquad(\bx\in\CY).
\end{displaymath}
See also Baldi et al.~\cite[Theorem~1.4]{bmp}. We place the exponent on the residual, as in~\cite{tran}. By compactness, it may be decreased to lie in $(0,1]$ after adjusting the constant. For the residuals below, we take $\varphi=\rho$ and $\psi=\dist(\cdot,\CX)$.

\begin{lemma}[~\cite{tran,extension}]\label{lem:atomic}
Let $\CX\subseteq\CY\subseteq\R^n$ be nonempty compact sets, and let $A\subseteq\N^n$ be finite. Let $\rho:\CY\to\R_+$ be continuous and satisfy $\dist(\bx,\CX)\leq c_{\Lo}\rho(\bx)^{\Lo}$ on $\CY$, where $0<\Lo\leq1$. Suppose a probability measure $\mu$ on $\CY$ satisfies
\begin{displaymath}
 \left\|\biy_A-\int_{\CY}\bv_A\,d\mu\right\|_\infty\leq e,
 \qquad \int_\CY\rho^2\,d\mu\leq\eta.
\end{displaymath}
If each coordinate of $\bv_A$ is $L_A$-Lipschitz on $\CY$, there is $\nu\in\PP(\CX)$ such that
\begin{equation}\label{eq:atomicrecall}
 \left\|\biy_A-\int_{\CX}\bv_A\,d\nu\right\|_\infty
 \leq e+L_Ac_{\Lo}\eta^{\Lo/2}.
\end{equation}
\end{lemma}
\begin{proof}
Consider the continuous map
\begin{displaymath}
 \bm{\Phi}:\CY\longrightarrow\R^{|A|+1},\qquad
 \bm{\Phi}(\bx)=\bigl(\bv_A(\bx),\rho(\bx)^2\bigr).
\end{displaymath}
Since $\CY$ is compact, the convex hull of $\bm{\Phi}(\CY)$ is compact and contains $\int_\CY\bm{\Phi}\,d\mu$. Carath\'eodory's theorem therefore gives points $\bx^{(1)},\ldots,\bx^{(s)}\in\CY$ and weights $\vartheta_b\ge0$, with $s\le|A|+2$ and $\sum_b\vartheta_b=1$, such that
\begin{equation}\label{eq:atomicpreservation}
 \sum_{b=1}^s\vartheta_b\bv_A(\bx^{(b)})=\int_\CY\bv_A\,d\mu,
 \qquad
 \sum_{b=1}^s\vartheta_b\rho(\bx^{(b)})^2=\int_\CY\rho^2\,d\mu\le\eta.
\end{equation}
This step preserves the residual integral even when $\rho^2$ is not a polynomial.

For each $b$, compactness of $\CX$ gives a nearest point $\bv^{(b)}\in\CX$ to $\bx^{(b)}$. Define the probability measure $\nu=\sum_b\vartheta_b\delta_{\bv^{(b)}}$. The assumed error bound and Jensen's inequality for the concave function $t\mapsto t^{\Lo/2}$ yield
\begin{displaymath}
\begin{aligned}
 \sum_b\vartheta_b\norm{\bx^{(b)}-\bv^{(b)}}_2
 &=\sum_b\vartheta_b\dist(\bx^{(b)},\CX)\\
 &\le c_{\Lo}\sum_b\vartheta_b\bigl(\rho(\bx^{(b)})^2\bigr)^{\Lo/2}
 \le c_{\Lo}\eta^{\Lo/2}.
\end{aligned}
\end{displaymath}
Using~\eqref{eq:atomicpreservation} and the coordinatewise Lipschitz bound, we obtain
\begin{displaymath}
 \left\|\int_\CY\bv_A\,d\mu-\int_\CX\bv_A\,d\nu\right\|_\infty
 \le L_A\sum_b\vartheta_b\norm{\bx^{(b)}-\bv^{(b)}}_2
 \le L_Ac_{\Lo}\eta^{\Lo/2}.
\end{displaymath}
The triangle inequality with the given approximation of $\biy_A$ proves~\eqref{eq:atomicrecall}.
\end{proof}

\paragraph{Choice of the residual function.}
After an affine change of variables, take $\CY=[0,1]^n$ and include the bounding-box generators in the hierarchy. The following choices specify $\rho$ in Lemma~\ref{lem:atomic} for all applications in this paper.

\emph{The hypercube.} When $\CX=\CY=[0,1]^n$, choose $\rho\equiv0$. Then $\eta=0$, no projection is needed, and the moment approximation error is $e$.

\emph{Scalar polynomial constraints.} For the description~\eqref{eq:scalarset}, choose
\begin{equation}\label{eq:scalarresidualrecall}
 \rho(\bx)=\left(\sum_{j=1}^{m_g}\max\{0,-g_j(\bx)\}^2
                     +\sum_{\ell=1}^{m_h}h_\ell(\bx)^2\right)^{1/2}.
\end{equation}
The functions $\rho$ and $\dist(\cdot,\CX)$ are continuous and semialgebraic on $\CY$, and their zero sets there both equal $\CX$. The recalled \L{}ojasiewicz inequality therefore gives
\begin{displaymath}
 \dist(\bx,\CX)\le c_{\Lo}\rho(\bx)^{\Lo}
 \qquad(\bx\in\CY),\qquad c_{\Lo}>0,\quad0<\Lo\le1.
\end{displaymath}
The equality-only case is obtained by setting $m_g=0$. Section~\ref{sec:projection} writes the same residual with clique indices and estimates its square by slack variable lifting.

\emph{Polynomial matrix constraints.} For the description~\eqref{eq:matrixPOP}, choose
\begin{equation}\label{eq:matrixresidualrecall}
 \rho_G(\bx)=\norm{G(\bx)_-}_{\mathrm F}.
\end{equation}
For a symmetric matrix $B$, $B_-$ has the same eigenvectors as $B$ and eigenvalues $\max\{0,-\lambda_i(B)\}$. Consequently
\begin{displaymath}
 \norm{B_-}_{\mathrm F}^2=\sum_i\max\{0,-\lambda_i(B)\}^2
 =\dist_{\mathrm F}(B,\mathbb S_+^p)^2
 \qquad(B\in\mathbb S^p).
\end{displaymath}
Use $\rho=\rho_G$ in Lemma~\ref{lem:atomic}. For the block-diagonal description $G=\diag(G_{aj})$, its square is the sum of the squared Frobenius norms of the negative parts of the blocks. The kernel variance estimate in Section~\ref{sec:matrix} controls the integral of this residual directly.

In the matrix setting, Dinh and Pham~\cite[Theorem~4.1]{dinhpham} give a \L{}ojasiewicz inequality in terms of the smallest eigenvalue, recalled by Tran and Toh~\cite[Theorem~2.4]{tranmatrix}. For the nonempty set $\CX=\{\bx:G(\bx)\succeq0\}$ and a compact ambient set $\CY$, there exist $c_{\Lo}>0$ and $0<\Lo\le1$ such that
\begin{displaymath}
 \dist(\bx,\CX)\le c_{\Lo}\max\{0,-\lambda_{\min}(G(\bx))\}^{\Lo}
 \le c_{\Lo}\rho_G(\bx)^{\Lo}\qquad(\bx\in\CY).
\end{displaymath}
Section~\ref{sec:matrix} applies this result to the block description with additional equalities and box constraints. In both settings, the constants and exponent concern the whole constraint system on the fixed compact ambient set. For $A=A_k$, each monomial has Lipschitz constant at most $k$ on $[0,1]^n$, so Lemma~\ref{lem:atomic} applies with $L_A=k$.

\paragraph{Polynomial lifting for scalar constraints.}
For a polynomial inequality, the negative-part residual $\max\{0,-g_j\}$ is not a polynomial. The lifting in~\cite{extension} makes it possible to estimate this residual using polynomial equations. For the scalar description~\eqref{eq:scalarset}, after a positive scaling, assume $\mon{g_j}\leq1$ and take $\CY=[0,1]^n$. Introduce $u_j\in[0,1]$ and consider
\begin{displaymath}
 \widehat{\CX}=\{(\bx,\bu)\in[0,1]^{n+m_g}:g_j(\bx)-u_j=0\ (j\in[m_g]),\ h_\ell(\bx)=0\ (\ell\in[m_h])\}.
\end{displaymath}
Its projection onto the $\bx$ coordinates is $\CX$. For a truncated functional $\ell_{\biy}$, put $\bg=(g_1,\ldots,g_{m_g})$ and define the lifted functional by polynomial substitution:
\begin{equation}\label{eq:denseliftrecall}
 \ell_{\widehat{\biy}}(\bx^\alpha \bu^\beta)=\ell_{\biy}(\bx^\alpha \bg(\bx)^\beta).
\end{equation}
The order of $\widehat{\biy}$ is chosen so that every substituted polynomial lies within the degree range of $\biy$. Positivity on the lifted box also requires a certificate for $g_j(1-g_j)$; we give the certificate and its degree bound in Lemma~\ref{lem:lift}.

The lifted functional annihilates the squared equations $(g_j-u_j)^2$. It also annihilates $h_\ell^2$ when the original equality ideal is enforced to the required degree. Therefore, if its degree-$t$ moments are approximated by a probability measure $\Gamma$ on the lifted hypercube with error $e_t$, then
\begin{displaymath}
 \int\left(\sum_j(g_j(\bx)-u_j)^2+\sum_\ell h_\ell(\bx)^2\right)d\Gamma
 \leq C e_t.
\end{displaymath}
Since $u_j\geq0$, this controls the squared feasibility residual of the original $\bx$ coordinates. Lemma~\ref{lem:atomic} then gives an error of the form
\begin{displaymath}
e_k+C'e_t^{\Lo/2}.
\end{displaymath}
This explains why a base error $e_t=O(r^{-2})$ becomes $O(r^{-\Lo})$, and why a factor $(\log_2r)^3$ becomes $(\log_2r)^{3\Lo/2}$.

\section{Gluing lemmas and transfer of convergence rates}\label{sec:method}
Using the notation and clique tree fixed in Section~\ref{sec:setup}, we prove two general gluing estimates. The first controls clique marginal errors in total variation; exact gluing follows as a corollary. The second controls feasibility errors and projects the joint measure onto a compact set. After explaining the role of running intersection, we combine these estimates into one rate-transfer principle with applications presented in Sections~\ref{sec:box}--\ref{sec:matrix}.

\subsection{Quantitative tree gluing}
\begin{lemma}[Tree gluing with total-variation error]\label{lem:gluing}
Let $\II$ have the rooted clique tree fixed above, let $\mu_a\in\PP(\HH^{I_a})$, and put
\begin{displaymath}
 \gamma_e=\TV((\pi_{S_e})_\#\mu_a,(\pi_{S_e})_\#\mu_b),
 \qquad e=\{a,b\}\in\mathscr T.
\end{displaymath}
There is $\Gamma\in\PP([0,1]^n)$ such that
\begin{equation}\label{eq:gluing}
 \TV((\pi_{I_a})_\#\Gamma,\mu_a)
 \le\sum_{e\in\operatorname{path}(\text{root},a)}\gamma_e.
\end{equation}
\end{lemma}
\begin{proof}
Start with the prescribed measure on the root clique and attach the other cliques in parent-before-child order. Let $a$ be the next child, let $b$ be its parent, and let $U$ be the union of the already attached coordinates. Running intersection gives $I_a\cap U=S:=I_a\cap I_b$. Write $E=[0,1]^S$, $F=[0,1]^{I_a\setminus S}$, and $\mu_{a,S}=(\pi_S)_\#\mu_a$.

Disintegrate $\mu_a$ with respect to its separator marginal. Such a disintegration exists because $E$ and $F$ are compact metric spaces; see also~\cite[Theorems~3.4 and~8.5]{kallenberg} and Chang and Pollard~\cite[Definition~1 and Theorems~1--2]{changpollard}. Thus choose a probability kernel $P_a$ from $E$ to $F$: $P_a(\bx_S,\cdot)$ is a probability measure for each $\bx_S$, and $\bx_S\mapsto P_a(\bx_S,B)$ is Borel measurable for each Borel $B\subseteq F$. The disintegration identity is
\begin{equation}\label{eq:disintegration}
 \int_{E\times F}\varphi(\bx_S,\bx_{I_a\setminus S})\,d\mu_a
 =\int_E\int_F\varphi(\bx_S,\bx_{I_a\setminus S})P_a(\bx_S,d\bx_{I_a\setminus S})\,\mu_{a,S}(d\bx_S)
\end{equation}
for every bounded Borel $\varphi$. These kernels are also called Markov kernels; see Kallenberg~\cite{kallenberg}. 

If $\nu$ is the existing measure on $[0,1]^U$, define its extension by
\begin{displaymath}
 \widetilde\nu(d\bx_U,d\bx_{I_a\setminus S})=\nu(d\bx_U)P_a(\bx_S,d\bx_{I_a\setminus S}).
\end{displaymath}
The kernel has mass one, so the $U$-marginal remains $\nu$. Moreover, it covers the coordinates $U \cup I_a$. The new clique marginal is $\nu_S\otimes P_a$, where $\nu_S=(\pi_S)_\#\nu$ and the notation denotes the joint measure including the separator coordinate. In particular, $\mu_a=\mu_{a,S}\otimes P_a$ by~\eqref{eq:disintegration}.

For any probability measures $\rho,\sigma$ on $E$ and the same probability kernel $P_a$,
\begin{equation}\label{eq:kernelTV}
 \TV(\rho\otimes P_a,\sigma\otimes P_a)=\TV(\rho,\sigma).
\end{equation}
Indeed, for a Borel set $C\subseteq E\times F$, the function $\bx_S\mapsto\int_F\mathbf1_C(\bx_S,\bx_{I_a\setminus S})P_a(\bx_S,d\bx_{I_a\setminus S})$ is measurable and takes values in $[0,1]$. Integrating it against $\rho-\sigma$ gives the upper bound in~\eqref{eq:kernelTV}. Taking $C=A\times F$ gives the reverse bound. Consequently,
\begin{displaymath}
\begin{aligned}
 \TV((\pi_{I_a})_\#\widetilde\nu,\mu_a)
 &=\TV(\nu_S,\mu_{a,S})\\
 &\le\TV(\nu_{I_b},\mu_b)+\gamma_{\{a,b\}},
\end{aligned}
\end{displaymath}
where the last step uses contraction under coordinate projection. Induction gives the path sum in~\eqref{eq:gluing}, and later attachments preserve it. The construction takes place on full hypercubes; an approximate change of separator measure need not preserve arbitrary local feasibility sets.
\end{proof}

When the separator marginals coincide, the same construction preserves every clique marginal and its support.

\begin{corollary}[Exact gluing and preservation of support]\label{cor:exactcompact}
Let $K_a\subseteq[0,1]^{I_a}$ be compact, and let $\mu_a\in\PP(K_a)$. Suppose the variable subsets form a clique tree and the separator marginals agree on each edge. Then there is a probability measure $\nu$ supported on
\begin{displaymath}
 K=\{\bx\in[0,1]^n:\bx_{I_a}\in K_a\text{ for every }a\}
\end{displaymath}
with $(\pi_{I_a})_\#\nu=\mu_a$ for every $a$. In particular, $K$ is nonempty.
\end{corollary}
\begin{proof}
Regard each $\mu_a$ as a probability measure on its hypercube. Lemma~\ref{lem:gluing}, with $\gamma_e=0$, produces a joint probability measure with the prescribed marginals. For each $a$, the cylinder $\{\bx:\bx_{I_a}\in K_a\}$ has probability one. Their finite intersection is $K$, so the joint measure is supported on $K$.
\end{proof}

\subsection{The main gluing lemma with feasibility errors}
Approximate gluing may move a clique marginal outside its local feasible set. The next lemma bounds the resulting feasibility error and projects the joint measure onto the global compact set. It permits auxiliary coordinates, so the same statement applies to polynomial lifts. Its two inputs are local moment errors and local residual integrals.

\begin{lemma}\label{lem:compactgluing}
Let $J_1,\ldots,J_m$ form a clique tree on $[N]$, where $n\le N$. Write a lifted point as $\widehat{\bx}=(\bx,\bu)\in[0,1]^N$, with original variables $\bx=(x_1,\ldots,x_n)$ and auxiliary variables $\bu\in[0,1]^{N-n}$. Let $\CX\subseteq[0,1]^n$ be nonempty and compact. Suppose continuous functions $\Phi_a:[0,1]^{J_a}\to\R_+$ satisfy
\begin{equation}\label{eq:compacterrorbound}
 \dist(\bx,\CX)\le c_{\Lo}
       \left(\sum_a\Phi_a(\widehat{\bx}_{J_a})\right)^{\Lo/2}
 \quad(\widehat{\bx}\in[0,1]^N),\qquad 0<\Lo\le1.
\end{equation}
Let $\mu_a\in\PP([0,1]^{J_a})$, and define, for every edge $e=\{a,b\}$,
\begin{displaymath}
\begin{aligned}
 \gamma_e&=\TV((\pi_{J_a\cap J_b})_\#\mu_a,
                  (\pi_{J_a\cap J_b})_\#\mu_b),\\
 \delta_a&=\sum_{e\in\operatorname{path}(\mathrm{root},a)}\gamma_e,
 \qquad \delta_*=\max_a\delta_a,\\
 \eta_a&=\int\Phi_a\,d\mu_a,\qquad
 M_a=\sup\Phi_a-\inf\Phi_a,\qquad
 \eta=\sum_a(\eta_a+M_a\delta_a).
\end{aligned}
\end{displaymath}
Let $A\subseteq\N^n$ be finite, with $|\alpha|\le k$ for every $\alpha\in A$, and suppose each support $\supp\alpha$ is contained in at least one $J_a\cap[n]$. Let $\biy_A$ be a prescribed real array satisfying
\begin{displaymath}
 \left|y_\alpha-\int \bx^\alpha d\mu_a\right|\le a_k
 \quad\bigl(\alpha\in A,\ \supp\alpha\subseteq J_a\bigr).
\end{displaymath}
Then there is a probability measure $\nu$ supported on $\CX$ such that
\begin{equation}\label{eq:compactgluing}
 \left\|\biy_A-\left(\int_{\CX}\bx^\alpha d\nu\right)_{\alpha\in A}\right\|_\infty
 \le a_k+\delta_*+k c_{\Lo}\eta^{\Lo/2}.
\end{equation}
\end{lemma}
\begin{proof}
Lemma~\ref{lem:gluing} gives a joint probability measure $\Gamma$ on $[0,1]^N$ whose $J_a$-marginal differs from $\mu_a$ by at most $\delta_a$ in total variation. Therefore
\begin{displaymath}
 \left|y_\alpha-\int \bx^\alpha d\Gamma\right|\le a_k+\delta_*,
 \qquad
 \int\sum_a\Phi_a(\widehat{\bx}_{J_a})\,d\Gamma\le\eta.
\end{displaymath}
The first inequality follows from~\eqref{eq:momentmetrics}. For the second, use the bound
\begin{displaymath}
 \left|\int v\,d\mu-\int v\,d\nu\right|
 \le(\sup v-\inf v)\TV(\mu,\nu)
\end{displaymath}
for bounded real $v$, with the convention $\TV(\mu,\nu)=\sup_B|\mu(B)-\nu(B)|$. By the finite-dimensional convex-hull theorem, replace $\Gamma$ by finitely many atoms $\widehat{\bx}^{(b)}=(\bx^{(b)},\bu^{(b)})$ with weights $\vartheta_b$ preserving all the target moments and the integral of $\sum_a\Phi_a$. For each atom, choose $\bv^{(b)}\in \CX$ nearest to $\bx^{(b)}$. Equation~\eqref{eq:compacterrorbound} and concavity give
\begin{displaymath}
\begin{aligned}
 \sum_b\vartheta_b\|\bx^{(b)}-\bv^{(b)}\|_2
 &\le c_{\Lo}\sum_b\vartheta_b
       \left(\sum_a\Phi_a(\widehat{\bx}_{J_a}^{(b)})\right)^{\Lo/2}\\
 &\le c_{\Lo}\left(\sum_b\vartheta_b
                   \sum_a\Phi_a(\widehat{\bx}_{J_a}^{(b)})\right)^{\Lo/2}
 \le c_{\Lo}\eta^{\Lo/2}.
\end{aligned}
\end{displaymath}
Set $\nu=\sum_b\vartheta_b\delta_{\bv^{(b)}}$. Each degree-$k$ monomial is $k$-Lipschitz on $[0,1]^n$. Its initial moment error and its change under projection give~\eqref{eq:compactgluing}.
\end{proof}

To apply this lemma without auxiliary variables, one may take
\begin{displaymath}
 \Phi_a(\bx_{I_a})=
 \sum_j\max\{0,-g_{aj}(\bx_{I_a})\}^2+
 \sum_\ell h_{a\ell}(\bx_{I_a})^2.
\end{displaymath}
Then~\eqref{eq:compacterrorbound} is precisely the global feasibility error bound. If the local measures are supported on their local feasible sets, all $\eta_a$ vanish. The bound still includes the separator errors $\delta_a$, because approximate gluing need not preserve feasibility. For pseudo-moment sequences, these negative parts cannot be evaluated by the Riesz functional. The polynomial lift in Section~\ref{sec:projection} replaces them by squared polynomial equations.

Exact separator compatibility makes $\delta_a=0$, but does not by itself make the local residual integrals $\eta_a$ vanish. Thus the kernel measures used below may still require projection onto $\CX$, even when their marginals glue exactly.

\subsection{Why the overlap structure is needed}\label{sec:structure}
The running-intersection property is the structural step behind the measure gluing lemma. When clique $I_a$ is attached after its parent, it gives
\begin{displaymath}
 I_a\cap\bigcup_{b<a}I_b=I_a\cap I_{p(a)}.
\end{displaymath}
Consequently, every previously assigned coordinate in $I_a$ is contained in one separator whose measure is already available. A conditional extension can then preserve the old joint measure. Without this property, different overlaps with previous cliques can impose mutually inconsistent joint requirements, even when all pairwise marginal tests agree. In the following example, we show that such gluing can fail without the running-intersection property. 

\begin{example}\label{ex:cycle}
Consider the variable cliques $\{1,2\}$, $\{2,3\}$, $\{3,4\}$, and $\{4,1\}$. On each of the first three cliques, we define the probability measure
$\frac12\delta_{(0,0)}+\frac12\delta_{(1,1)}$, and on the last clique put
$\frac12\delta_{(0,1)}+\frac12\delta_{(1,0)}$. Every one-variable marginal is the uniform probability measure on $\{0,1\}$, so all overlap measures, and hence all overlap moments, agree exactly. Nevertheless, no joint probability measure has these marginals. Indeed, the first three marginals require
$x_1=x_2=x_3=x_4$ almost surely, while the last requires $x_4\ne x_1$ almost surely.

The four-cycle graph is not chordal, and its edge cliques have no running-intersection ordering. A chordal extension must enlarge the local cliques and enforce the resulting larger separator constraints. It cannot be used merely to relabel these four inconsistent measures. Thus running intersection cannot be omitted from the general gluing guarantee on the basis of pairwise consistency alone. This does not assert that a particular tuple on a cyclic pattern can never glue; independent product marginals, for example, do glue.
\end{example}

A chordal extension is a way to obtain maximal variable cliques with a running-intersection ordering. It may enlarge the blocks and the separators, and the relaxation must enforce the constraints associated with those enlarged sets. Once the clique tree has been fixed, the proof uses its running-intersection property directly. This structural requirement concerns measure gluing; it does not imply same-order completion of pseudo-moments.

\subsection{Transferring the dense approximation rate}
The gluing lemmas presented above play a role as tools to transfer the convergence rates and the Hausdorff distance from the dense cases to the sparse cases that we present in what follows. 
\begin{theorem}\label{thm:ratetransfer}
Fix a clique tree and a target degree $k$. Suppose that, for each sparse pseudo-moment sequence at order $r$, a local construction gives probability measures $\mu_a$ on the clique hypercubes with degree-$k$ moment error at most $a_k(r)$. Let $\gamma_e(r)$ be their separator discrepancies in total variation and put
\begin{displaymath}
\delta_*(r)=\max_a\sum_{e\in\operatorname{path}(\mathrm{root},a)}\gamma_e(r).
\end{displaymath}
Then one joint probability measure $\Gamma$ satisfies
\begin{equation}\label{eq:abstractambient}
 \max_{\alpha\in A_k}
    \left|y_\alpha-\int \bx^\alpha\,d\Gamma\right|
 \le a_k(r)+\delta_*(r).
\end{equation}
The construction may be applied in original or lifted coordinates, with $A_k$ indexing only original-coordinate monomials.

For a nonempty compact $\CX\subseteq[0,1]^n$, suppose additionally that a continuous residual $\rho\ge0$ satisfies
\begin{displaymath}
 \dist(\bx,\CX)\le c_{\Lo}\rho(\bx)^{\Lo},\qquad 0<\Lo\le1,
 \qquad \int\rho(\bx)^2\,d\Gamma\le\eta(r).
\end{displaymath}
There is $\nu\in\PP(\CX)$ with truncated-moment error at most
\begin{equation}\label{eq:abstractcompact}
 a_k(r)+\delta_*(r)+k c_{\Lo}\eta(r)^{\Lo/2}.
\end{equation}
If the estimates are uniform over sparse pseudo-moment sequences and exact moments are feasible, these are Hausdorff bounds. In particular, $a_k,\delta_*=O(\psi(r))$ give the hypercube rate $O(\psi(r))$, and $a_k,\delta_*,\eta=O(\psi(r))$ give the compact-set rate $O(\psi(r)^{\Lo/2})$, where $\psi(r)\to0$.
\end{theorem}
\begin{proof}
Apply Lemma~\ref{lem:gluing}. Equation~\eqref{eq:momentmetrics} and the triangle inequality give~\eqref{eq:abstractambient}; the same joint measure works for all target degrees covered by the local estimates. The atomic projection estimate, Lemma~\ref{lem:atomic}, applied to its $\bx$-marginal gives~\eqref{eq:abstractcompact}. Lemma~\ref{lem:compactgluing} supplies the required residual bound from local residual integrals and separator discrepancies. Taking the supremum over sparse pseudo-moment sequences gives what we desire. 
\end{proof}

In the following sections, we apply the gluing lemma that we have developed in this section to present the method to convey the convergence rates and bounds on the Hausdorff distances from the dense cases to the sparse cases. In each section, we also recall the existing results in the mentioned dense cases for clarity. 

\section{Convergence rates on the hypercube}\label{sec:box}
We first treat the correlative sparse preordering. Its normalized Jackson kernel preserves separator marginals exactly, so a short kernel-gluing lemma transfers the dense approximation rate without an additional error. We then treat the correlative sparse quadratic module: the squared kernel is only approximately normalized, and we control the resulting separator discrepancies in total variation before applying the measure-gluing lemma. The kernel estimates themselves are recalled from Laurent and Slot~\cite{laurent} and Gribling et al.~\cite{kernel}.

\subsection{The correlative sparse preordering}\label{sec:boxT}
We first recall the Jackson kernel, its degree bounds, and the positivity and approximation properties needed below. We then construct compatible local measures and derive the Hausdorff estimate by exact gluing.

\paragraph{The Jackson kernel and its properties.}
Use the Chebyshev probability measure $\varpi$, its products $\varpi_I$, and the shifted basis $\tau_j(x)=T_j(2x-1)$ from Section~\ref{sec:momentsets}. Laurent and Slot~\cite[Proposition~6, Lemma~7, and Theorem~8]{laurent} provide the following properties of the modified Jackson kernel, written here on $[0,1]$. For every even $q\ge2$, there is a polynomial kernel $J_q(t,x)$ of degree at most $q$ in each variable such that
\begin{equation}\label{eq:Jackson}
 \int J_q(t,x)\,d\varpi(t)=1,\qquad
 J_q(t,\cdot)\in\QQ_{\{i\},q}([0,1])\quad(t\in[0,1]).
\end{equation}
For a coordinate set $I$, use the same univariate kernel in
\begin{displaymath}
 J_I(\bt_I,\bx_I)=\prod_{i\in I}J_q(t_i,x_i),\qquad
 \mathcal J_Iv(\bx_I)=\int v(\bt_I)J_I(\bt_I,\bx_I)\,d\varpi_I(\bt_I).
\end{displaymath}
The product kernel satisfies $J_I(\bt_I,\cdot)\in\TT_{I,|I|q}([0,1]^I)$. The coefficient and tensor estimates in~\cite[Proposition~6(iii) and Lemma~7]{laurent} give
\begin{equation}\label{eq:Jacksontensor}
\begin{gathered}
 \cheb{\mathcal J_Iv-v}\le\Delta_s(q)\cheb v,\qquad
 \Delta_s(q)=\frac{\pi^2s^2}{(q+2)^2},\\
 \deg(\mathcal J_Iv-v)\le s\qquad(\deg v\le s\le q).
\end{gathered}
\end{equation}
These normalization, certificate, degree, and approximation properties are the only kernel facts needed below.

\paragraph{Exact gluing and the Hausdorff bound.}
We now apply the recalled kernel to a sparse pseudo-moment sequence. Its preordering certificate gives a nonnegative density under each clique pseudo-moment functional, and its exact normalization makes the separator marginals coincide. The following lemma glues the resulting probability measures to approximate truncated sparse pseudo-moment sequences. 

\begin{lemma}\label{lem:compatiblekernel}
Let $\biy\in\FF_r^{\TT}([0,1]^n)$ be a sparse pseudo-moment sequence, with clique restrictions $\biy^{(a)}$. Use the same Jackson kernel $J_q$ in every coordinate and clique, with $q\ge2$ even and $wq\le2r$. Then
\begin{displaymath}
 d\mu_a(\bt_{I_a})=\ell_{\biy^{(a)}}(J_{I_a}(\bt_{I_a},\cdot))d\varpi_{I_a}(\bt_{I_a})
\end{displaymath}
are probability measures with identical separator marginals. They admit a joint probability measure $\Gamma$ whose $I_a$-marginal is $\mu_a$ for every $a$. Every simultaneous bound on the local moment-approximation errors therefore holds for $\Gamma$ with the same error.
\end{lemma}
\begin{proof}
Fix a clique $I_a$ and $\bt=\bt_{I_a}\in[0,1]^{I_a}$. The recalled preordering certificate gives
\begin{displaymath}
 J_{I_a}(\bt,\bx)=\sum_{A\subseteq I_a}\sum_\nu p_{A,\nu,\bt}(\bx)^2b_A(\bx),
 \qquad b_A(\bx)=\prod_{i\in A}x_i(1-x_i),
\end{displaymath}
where each summand has degree at most $|I_a|q\le2r$; in particular, $\deg p_{A,\nu,\bt}\le r-|A|$. Feasibility imposes
\begin{displaymath}
 M_{r-|A|}(b_A\biy^{(a)})\succeq0
 \qquad(A\subseteq I_a).
\end{displaymath}
If $\bc_{A,\nu,\bt}$ denotes the coefficient vector of $p_{A,\nu,\bt}$ in the monomial basis through degree $r-|A|$, then the density satisfies
\begin{displaymath}
 h_a(\bt):=\ell_{\biy^{(a)}}(J_{I_a}(\bt,\cdot))
 =\sum_{A,\nu}\bc_{A,\nu,\bt}^{T}
 M_{r-|A|}(b_A\biy^{(a)})\bc_{A,\nu,\bt}\ge0.
\end{displaymath}
The term $A=\varnothing$ uses the ordinary moment matrix. Thus positivity follows from the moment and preordering localizing matrices; pointwise nonnegativity of the kernel alone would not suffice for a pseudo-moment functional.

Since $h_a$ is a polynomial in $\bt$, it is measurable and bounded on the cube. Consequently $\mu_a(B)=\int_B h_a\,d\varpi_{I_a}$ defines a finite positive Borel measure. Exact kernel normalization and linearity give
\begin{displaymath}
 \mu_a([0,1]^{I_a})
 =\ell_{\biy^{(a)}}\left(\int J_{I_a}(\bt,\cdot)\,d\varpi_{I_a}(\bt)\right)
 =\ell_{\biy^{(a)}}(1)=1,
\end{displaymath}
so $\mu_a$ is a probability measure. For a separator $S=I_a\cap I_b$, integrating out the other coordinates leaves exactly $J_S$. Hence
\begin{displaymath}
 \frac{d(\pi_S)_\#\mu_a}{d\varpi_S}
 =\ell_{\biy^{(a)}}(J_S(\bt_S,\cdot))
 =\ell_{\biy^{(b)}}(J_S(\bt_S,\cdot))
 =\frac{d(\pi_S)_\#\mu_b}{d\varpi_S},
\end{displaymath}
where the middle equality uses separator consistency. Exact tree gluing, Corollary~\ref{cor:exactcompact}, gives $\Gamma$. Its prescribed clique marginals give the moment assertion.
\end{proof}

\begin{corollary}\label{thm:boxT}
For the correlative sparse box preordering, let $q\ge2$ be even, $wq\le2r$, and $k\le\min\{q,r\}$. Then
\begin{displaymath}
 \bd_{k,r}^{\TT}([0,1]^n)\le\Delta_k(q),\qquad
 0\le f_{\min}-p_r^{\TT}\le F_1\Delta_k(q).
\end{displaymath}
For every sparse pseudo-moment sequence, one probability measure realizes the corresponding moment bounds simultaneously for $1\le s\le\min(q,r)$. In particular, $q=2\lfloor r/w\rfloor$ gives $O(r^{-2})$ for fixed $k,w$.
\end{corollary}
\begin{proof}
For each sparse pseudo-moment sequence, take the local probability measures and their joint measure from Lemma~\ref{lem:compatiblekernel}. For $|\alpha|\le s\le\min(q,r)$,~\eqref{eq:Jacksontensor} and Lemma~\ref{lem:evaluation} give
\begin{displaymath}
 \left|\int \bt^\alpha d\mu_a-y_\alpha\right|
 =\left|\ell_{\biy^{(a)}}(\mathcal J_{I_a}\bx^\alpha-\bx^\alpha)\right|
 \le\Delta_s(q),
\end{displaymath}
since $\cheb{\bx^\alpha}=1$. Exact gluing preserves these bounds. Thus
$\biz=(\int \bx^\alpha d\Gamma)_{\alpha\in A_k}\in\MM_k^{\II}([0,1]^n)$ satisfies $\norm{\biy_{A_k}-\biz}_\infty\le\Delta_k(q)$ for the restriction $\biy_{A_k}$ of every sparse pseudo-moment sequence $\biy$. Taking the supremum in~\eqref{eq:distance} proves the Hausdorff bound, and the bound on the gap between $f_{\min}$ and $p_r^{\TT}$ follows. 
\end{proof}

\subsection{The correlative sparse quadratic module}\label{sec:boxQ}
We first recall the squared kernel and its SOS, mass, and approximation estimates. We then bound the separator discrepancies of its local probability measures in total variation and use measure gluing to derive the Hausdorff bound.

\paragraph{The squared kernel and its properties.}
We use the approximately normalized squared Fej\'er kernel of Gribling et al.~\cite[Sections~3.2 and~3.4]{kernel}, shifted to $[0,1]$. For integers $q\ge2$ and even $N\ge2$, put
\begin{equation}\label{eq:kernel}
 C_q=\frac{2q^2+1}{3q},\qquad D=2q(N+1),\qquad\delta=2^{-N-1}.
\end{equation}
Their kernel $\kappa(t,x)$ and its operator
\begin{displaymath}
 \mathcal Kv(x)=\int v(t)\kappa(t,x)\,d\varpi(t),\qquad m=\mathcal K1,
\end{displaymath}
have the following properties. For fixed $t$, $\kappa(t,\cdot)$ and $m$ are SOS polynomials of degree at most $D$. The mass estimates in~\cite[Lemmas~1--2, Section~3.2, and Proposition~8]{kernel}, together with the interval representation in~\cite[Theorem~8]{laurent}, give
\begin{equation}\label{eq:masscertificates}
 1-\delta\le m\le1\text{ on }[0,1],\qquad
 1-m,\ m-(1-\delta)\in\QQ_{\{i\},D}([0,1]).
\end{equation}
The univariate approximation estimates in~\cite[Lemma~5 and Propositions~7--8]{kernel} give the following bound, with conservative constants:
\begin{equation}\label{eq:eps}
\begin{gathered}
 \cheb{\mathcal K\tau_j-\tau_j}\le\varepsilon_s\qquad(0\le j\le s\le q),\\
 \varepsilon_s:=\frac{N+1}{C_q}\left(\frac{s^2}{q^2}+\frac{3s^2}{2q}\right)
          +2\sqrt{2(D+1)}\,\delta.
\end{gathered}
\end{equation}
For the product kernel and operator
\begin{displaymath}
 K_I(\bt_I,\bx_I)=\prod_{i\in I}\kappa(t_i,x_i),\qquad
 \mathcal K_Iv(\bx_I)=\int v(\bt_I)K_I(\bt_I,\bx_I)\,d\varpi_I(\bt_I),
\end{displaymath}
the tensor estimate~\cite[Lemma~10]{kernel} reads
\begin{equation}\label{eq:squaredtensor}
\begin{gathered}
 \cheb{\mathcal K_Iv-v}\le\bigl((1+\varepsilon_s)^{|I|}-1\bigr)\cheb v,\\
 \deg(\mathcal K_Iv-v)\le|I|D\qquad(\deg v\le s\le q).
\end{gathered}
\end{equation}
We use the same $q,N$ and univariate kernel in every coordinate and clique.

\paragraph{Why the squared kernel requires gluing with error.}
For the modified Jackson kernel in Lemma~\ref{lem:compatiblekernel}, integrating out a coordinate gives the constant polynomial $1$. Thus integrating out all coordinates outside a separator $S=I_a\cap I_b$ leaves exactly $J_S$, and separator moment consistency makes the two marginal densities identical. Its tensor-product positivity certificate, however, lies in the box preordering; it need not belong to the box quadratic module, which does not contain arbitrary products of its generators.

The squared kernel has the SOS certificate needed for the quadratic-module hierarchy, but its mass is the polynomial $m(x)$ rather than the constant $1$. Put $M_I(\bx_I)=\prod_{i\in I}m(x_i)$. Marginalization now gives
\begin{displaymath}
 \int_{[0,1]^{I_a\setminus S}}K_{I_a}(\bt_{I_a},\bx_{I_a})
 \,d\varpi_{I_a\setminus S}(\bt_{I_a\setminus S})
 =K_S(\bt_S,\bx_S)M_{I_a\setminus S}(\bx_{I_a\setminus S}).
\end{displaymath}
The extra factor depends on variables outside $S$. Consistency identifies the two clique functionals only on polynomials supported in $S$, so it does not equate their evaluations of these products. Dividing each clique density by its scalar total mass $Z_a=\ell_{\biy^{(a)}}(M_{I_a})$ ensures probability mass one but does not cancel this variable-dependent factor. Therefore exact separator compatibility is not guaranteed.

Dividing the kernel pointwise by $m(x)$ would restore exact normalization, but in general produces a rational function, outside the finite-degree polynomial tests of the pseudo-moment functional. We instead use $1-\delta\le m\le1$ together with its quadratic-module certificates to bound the separator discrepancy. The next lemma gives the error $2w\delta$, which the quantitative gluing lemma accumulates along tree paths. This explains why this squared-kernel construction requires gluing with error; it does not rule out other kernels with exact compatibility.

\begin{lemma}\label{lem:squaredseparator}
Let $\biy\in\FF_r^{\QQ}([0,1]^n)$ be a sparse pseudo-moment sequence, with clique restrictions $\biy^{(a)}$, and suppose $wD\le2r$. Put $M_I=\prod_{i\in I}m(x_i)$ and $Z_a=\ell_{\biy^{(a)}}(M_{I_a})$. Then
\begin{displaymath}
 (1-\delta)^{|I_a|}\le Z_a\le1,\qquad
 d\mu_a(\bt)=Z_a^{-1}\ell_{\biy^{(a)}}(K_{I_a}(\bt,\cdot))d\varpi_{I_a}(\bt)
\end{displaymath}
defines a probability measure. On every edge $\{a,b\}$ with separator $S=I_a\cap I_b$,
\begin{equation}\label{eq:squaredseparator}
 \TV\bigl((\pi_S)_\#\mu_a,(\pi_S)_\#\mu_b\bigr)\le2w\delta.
\end{equation}
\end{lemma}
\begin{proof}
For disjoint $B,C,\{j\}\subseteq I_a$, the polynomial $U=K_BM_C$ is SOS. Multiply the two certificates in~\eqref{eq:masscertificates} by $U$. Their degrees are at most $wD\le2r$, so feasibility gives
\begin{equation}\label{eq:massinsert}
 (1-\delta)\ell_{\biy^{(a)}}(U)
 \le\ell_{\biy^{(a)}}(Um(x_j))\le\ell_{\biy^{(a)}}(U).
\end{equation}
Starting from $U=1$ gives the bounds on $Z_a$. Positivity of the SOS kernel then proves that $\mu_a$ is a probability measure.

Consistency gives a common nonnegative separator density
\begin{displaymath}
 h_S(\bt_S)=\ell_{\biy^{(a)}}(K_S(\bt_S,\cdot))
         =\ell_{\biy^{(b)}}(K_S(\bt_S,\cdot)),\qquad
 z_S=\int h_Sd\varpi_S\ge(1-\delta)^{|S|}>0.
\end{displaymath}
The actual separator marginal is
\begin{equation}\label{eq:marginaldefect}
 \frac{d(\pi_S)_\#\mu_a}{d\varpi_S}
 =Z_a^{-1}\ell_{\biy^{(a)}}(K_SM_{I_a\setminus S}).
\end{equation}
Repeated use of~\eqref{eq:massinsert} bounds the numerator between $(1-\delta)^{|I_a\setminus S|}h_S$ and $h_S$, and gives $Z_a\le z_S$. Thus
\begin{displaymath}
 (\pi_S)_\#\mu_a\ge(1-\delta)^{|I_a\setminus S|}\nu_S,
 \qquad d\nu_S=(h_S/z_S)d\varpi_S.
\end{displaymath}
For probability measures $\mu\ge c\nu$, $\TV(\mu,\nu)\le1-c$. Each marginal is therefore within $w\delta$ of $\nu_S$, then~\eqref{eq:squaredseparator} follows from the triangle inequality.
\end{proof}

Lemma~\ref{lem:squaredseparator} controls the normalization defect directly at the level of measures; it does not infer total variation from finitely many moment differences.

\begin{corollary}\label{thm:boxQ}
For a sparse pseudo-moment sequence of the correlative sparse box quadratic-module hierarchy, choose $q\ge\max\{2,k\}$ and even $N\ge2$ with
\begin{equation}\label{eq:budget}
 wD=2wq(N+1)\le r.
\end{equation}
There is one probability measure $\Gamma$ on $[0,1]^n$ satisfying, simultaneously for $1\le s\le q$,
\begin{equation}\label{eq:E}
 \max_{\alpha\in A_s}\left|y_\alpha-\int \bx^\alpha d\Gamma\right|
 \le E_s(q,N;w,\omega):=\left(\frac{1+\varepsilon_s}{1-\delta}\right)^w-1+2w\omega\delta.
\end{equation}
Consequently $\bd_{k,r}^{\QQ}([0,1]^n)\le E_k$ and
$0\le f_{\min}-p_r^{\QQ}\le F_1E_k$. Both are $O((\log_2r)^3/r^2)$ for fixed data.
\end{corollary}
\begin{proof}
Take the local probability measures from Lemma~\ref{lem:squaredseparator}. By~\eqref{eq:squaredtensor} and Lemma~\ref{lem:evaluation}, the unnormalized moment error is at most $\beta_s=(1+\varepsilon_s)^w-1$: its polynomial degree is at most $wD\le r$ and $\cheb{\bx^\alpha}=1$. Since $|y_\alpha|\le1$,
\begin{displaymath}
 \left|\int \bt^\alpha d\mu_a-y_\alpha\right|
 \le\frac{\beta_s+1-Z_a}{Z_a}
 \le\left(\frac{1+\varepsilon_s}{1-\delta}\right)^w-1.
\end{displaymath}
Apply Lemma~\ref{lem:gluing} with separator errors $2w\delta$. It produces $\Gamma$ whose clique marginals differ from $\mu_a$ by at most $2w\omega\delta$ in total variation. Since $0\le \bt^\alpha\le1$, this adds at most $2w\omega\delta$ to each moment error and proves~\eqref{eq:E}. The construction applies to every sparse pseudo-moment sequence, so taking the supremum in~\eqref{eq:distance} gives $\bd_{k,r}^{\QQ}([0,1]^n)\le E_k$. This is the Hausdorff distance because true moments are feasible. Integrating the objective gives its gap bound.

Choose $N=2\lceil2\log_2(q+1)\rceil$ and the largest admissible $q$. Then $q\asymp r/(w\log_2r)$, $\delta\le\tfrac12(q+1)^{-4}$, and $\varepsilon_s=O_s((\log_2q)/q^2)$, proving the rate. The sufficient budget~\eqref{eq:budget} includes the normalization polynomial and the degree restriction of Lemma~\ref{lem:evaluation}; the hierarchy itself uses certificate degree $2r$.
\end{proof}

\section{Convergence rates on compact sets}\label{sec:projection}
We now consider the scalar description~\eqref{eq:scalarset} and combine the hypercube rates with Lemma~\ref{lem:compactgluing} to approximate truncated moments on a compact basic semialgebraic set. The proof has three steps: lift each clique pseudo-moment restriction using slack variables, apply the hypercube kernels, and glue and project the resulting local measures. This section verifies the degree and residual estimates needed for that application.

\subsection{Polynomial lifting and the degree bound}

We now apply the method recalled in Section~\ref{sec:recall} to a sparse pseudo-moment sequence. The substitution and the atomic projection are the same as in~\cite{tran,extension}; we verify their degree requirements and their compatibility with the clique tree.
Scale each nonzero $g_{aj}$ by a positive constant to have $\mon{g_{aj}}\le1$. Put $d=\max(\{1\}\cup\{\deg g_{aj}:a,j\})$. Give every inequality a slack v a variable $u_{aj}$, and set
\begin{displaymath}
 \widehat I_a=I_a\cup U_a,\qquad\widehat w=\max_a|\widehat I_a|.
\end{displaymath}
All slacks are distinct across cliques. The separators and tree depth remain unchanged. Write $\bg_a=(g_{aj})_j$. For a polynomial on a lifted clique, define
\begin{equation}\label{eq:lift}
 \ell_{\widehat{\biy}^{(a)}}(p)=
 \ell_{\biy^{(a)}}\bigl(p(\bx_{I_a},\bg_a(\bx_{I_a}))\bigr).
\end{equation}

\begin{lemma}[Positivity of the polynomial lift]\label{lem:lift}
Let $\biy\in\FF_r^{\mathcal C}(\CX)$ be a sparse pseudo-moment sequence. Its clique lifts in~\eqref{eq:lift} define a sparse pseudo-moment sequence $\widehat{\biy}$ on the lifted hypercube at order $R$ if
\begin{displaymath}
 d(R+1)\le r\quad(\QQ),\qquad 2dR\le r\quad(\TT).
\end{displaymath}
It preserves normalization and full separator consistency.
\end{lemma}
\begin{proof}
Substitution of degree $2R$ has degree at most $2dR$. Lemma~\ref{lem:evaluation}'s certificate proof gives
\begin{displaymath}
 1-g\in\QQ_{I_a,2\lceil\deg g/2\rceil}([0,1]^{I_a})\subseteq\QQ_{I_a,2d}.
\end{displaymath}
The crucial certificate is
\begin{equation}\label{eq:slack}
 g(1-g)=g^2(1-g)+(1-g)^2g\in\QQ_{a,4d}(\CX).
\end{equation}
Only SOS multipliers have been used. A lifted SOS substitutes to an SOS. A slack-box localizer of degree at most $2R$ has, after substitution and \eqref{eq:slack}, certificate degree at most $4d+d(2R-2)=2dR+2d\le2r$. Original box localizers have no larger budget. This proves the quadratic-module assertion.

For a local preordering term $\sigma\prod_{i\in A}b_i\prod_{j\in B}u_j(1-u_j)$, the substituted certificate has degree at most
\begin{displaymath}
 d\deg\sigma+2|A|+4d|B|\le4dR\le2r.
\end{displaymath}
Each local preordering is closed under these products, with repeated generator factors absorbed into squares. Slacks variable do not occur in a separator, so \eqref{eq:consistency} is preserved.
\end{proof}

Set
\begin{align}
 t&=\max\{k,2d,2\deg h_{a\ell}:a,\ell\},\nonumber\\
 H_a(\bx,\bu)&=\sum_j(g_{aj}(\bx)-u_{aj})^2+\sum_\ell h_{a\ell}(\bx)^2,
 \qquad C_H=\sum_a\mon{H_a},\nonumber\\
 \rho(\bx)^2&=\sum_{a,j}\max\{0,-g_{aj}(\bx)\}^2+\sum_{a,\ell}h_{a\ell}(\bx)^2.\label{eq:residual}
\end{align}
The residual in~\eqref{eq:residual} is~\eqref{eq:scalarresidualrecall} written with the assigned clique indices. Apply the semialgebraic \L{}ojasiewicz inequality of Bochnak et al.~\cite[Corollary~2.6.7]{bcr}, recalled in Section~\ref{sec:recall}, to $\rho$ and $\dist(\cdot,\CX)$ on $[0,1]^n$. It provides a constant $c_{\Lo}>0$ and an exponent $0<\Lo\le1$ satisfying~\eqref{eq:LoIntro}, as in the scalar moment-approximation argument of Tran and Toh~\cite{tran}. Fix such a pair for the remainder of this section. Replacing this Euclidean residual by a sum or a maximum of the same violations changes constants only. The inequality holds for the whole intersection $\CX$ on the ambient hypercube; separate inequalities for the local feasible sets do not imply it with the same exponent.

The following theorem shows that the slack variables can be used to convey the above method to compact algebraic sets. 
\begin{theorem}\label{thm:projection}
Let $\biy\in\FF_r^{\mathcal C}(\CX)$, and let $R$ satisfy Lemma~\ref{lem:lift}, with $t\le2R$ and $2\deg h_{a\ell}\le2r$. Apply the common-kernel construction to the lifted sparse pseudo-moment sequence. Suppose the local probability measures have moment errors at most $a_k,a_t$, and let $\delta_*$ be the maximum path sum of their separator discrepancies. Put
\begin{displaymath}
e_k=a_k+\delta_*,\qquad e_t=a_t+\delta_*.
\end{displaymath}
Then there is $\nu\in\PP(\CX)$ such that
\begin{equation}\label{eq:maintransfer}
 \max_{\alpha\in A_k}\left|y_\alpha-\int_{\CX}\bx^\alpha d\nu\right|
 \le e_k+k c_{\Lo}(C_He_t)^{\Lo/2}.
\end{equation}
If $a_k,a_t$ and the bound on $\delta_*$ are uniform over $\biy\in\FF_r^{\mathcal C}(\CX)$, then the right side bounds $\bd_{k,r}^{\mathcal C}(\CX)$. Under the same uniformity assumption, for $f=\sum_a f_a(\bx_{I_a})$ with $\deg f_a\le k$,
\begin{equation}\label{eq:objective}
 0\le f_{\min,\CX}-p_r^{\mathcal C}(f;\CX)
 \le F_1e_k+\Lip(f)c_{\Lo}(C_He_t)^{\Lo/2}.
\end{equation}
Here $\Lip(f)$ is taken on $[0,1]^n$ in Euclidean distance.
\end{theorem}
\begin{proof}
Substitution gives
\begin{displaymath}
 \ell_{\widehat{\biy}^{(a)}}((g_{aj}-u_{aj})^2)=0,
 \qquad
 \ell_{\widehat{\biy}^{(a)}}(h_{a\ell}^2)=0.
\end{displaymath}
The first identity is algebraic. The second follows from the truncated equality ideal. Hence $\ell_{\widehat{\biy}^{(a)}}(H_a)=0$. Since $\deg H_a\le t$, the local moment estimate implies
\begin{displaymath}
 0\le\int H_a\,d\mu_a\le\mon{H_a}a_t.
\end{displaymath}
On the lifted hypercube $H_a\ge0$ and $\sup H_a\le\mon{H_a}$, so its oscillation is at most $\mon{H_a}$. Also $u_{aj}\ge0$ gives
\begin{displaymath}
 \rho(\bx)^2\le\sum_a H_a(\bx,\bu),\qquad
 \dist(\bx,\CX)\le c_{\Lo}\left(\sum_aH_a(\bx,\bu)\right)^{\Lo/2}.
\end{displaymath}
Apply Lemma~\ref{lem:compactgluing} to $J_a=\widehat I_a$ and $\Phi_a=H_a$. Its total feasibility bound satisfies
\begin{equation}\label{eq:Hmean}
 \eta\le\sum_a\mon{H_a}(a_t+\delta_a)
       \le C_H(a_t+\delta_*)=C_He_t.
\end{equation}
This proves~\eqref{eq:maintransfer}.

For the objective estimate, preserve the degree-$k$ target moments in the atomic representation used by the lemma. The initial error in the integral of $f$ is at most $F_1e_k$. Projection changes it by at most $\Lip(f)$ times the mean projection distance, which is bounded by $c_{\Lo}(C_He_t)^{\Lo/2}$. For this $\biy$, the inequality $\int_{\CX}f\,d\nu\ge f_{\min,\CX}$ gives
\begin{displaymath}
 f_{\min,\CX}-\ell_{\biy}(f)
 \le F_1e_k+\Lip(f)c_{\Lo}(C_He_t)^{\Lo/2}.
\end{displaymath}
When the error bounds are uniform over $\FF_r^{\mathcal C}(\CX)$, taking the infimum of $\ell_{\biy}(f)$ proves~\eqref{eq:objective}. Exact sparse moment sequences give the lower bound zero. The kernel estimates below provide precisely this uniformity, and none of the constants depends on the number of atoms.
\end{proof}

For the quadratic-module kernel, the proof of Corollary~\ref{thm:boxQ} gives
\begin{displaymath}
 a_s=\left(\frac{1+\varepsilon_s}{1-\delta}\right)^{\widehat w}-1,
 \qquad \delta_*\le2\widehat w\omega\delta.
\end{displaymath}
Thus $e_s$ can be replaced by $E_s(q,N;\widehat w,\omega)$ for $s=k,t$. For the preordering kernel, $a_s=\Delta_s(q)$ and $\delta_*=0$. This identifies separately the moment-approximation error, the gluing error, and the projection error in~\eqref{eq:maintransfer}.

\subsection{Convergence rates for correlative sparsity}

\begin{theorem}[Constrained correlative sparse rates]\label{thm:CS}
Under the stated box-generator and clique-tree assumptions, with $\Lo$ chosen from the scalar \L{}ojasiewicz inequality above, for fixed $k$ and fixed problem data,
\begin{displaymath}
\begin{aligned}
 \bd_{k,r}^{\QQ}(\CX),\ f_{\min,\CX}-p_r^{\QQ}(f;\CX)
 &=O\left((\log_2r)^{3\Lo/2}/r^{\Lo}\right),\\
 \bd_{k,r}^{\TT}(\CX),\ f_{\min,\CX}-p_r^{\TT}(f;\CX)
 &=O(r^{-\Lo}).
\end{aligned}
\end{displaymath}
The moment values and matching SOS suprema converge increasingly to $f_{\min,\CX}$.
\end{theorem}
\begin{proof}
For $\QQ$, take $R=\lfloor r/d\rfloor-1$. Apply Corollary~\ref{thm:boxQ} to the lifted clique tree with width $\widehat w$ and a common $(q,N)$ satisfying $q\ge t$ and $2\widehat wq(N+1)\le R$. Then $e_s=E_s(q,N;\widehat w,\omega)$ for $s=k,t$. The parameter choice in its proof gives $e_k,e_t=O((\log_2r)^3/r^2)$ since $R\asymp r$. Insert this into \eqref{eq:maintransfer}--\eqref{eq:objective}.

For $\TT$, take $R=\lfloor r/(2d)\rfloor$ and $q=2\lfloor R/\widehat w\rfloor$. For large $r$, $t\le\min(q,R)$, and Corollary~\ref{thm:boxT} gives $e_s=\Delta_s(q)=O(r^{-2})$. The projection exponent is again $\Lo/2$. Restriction from order $r+1$ to order $r$ is feasible, proving monotonicity.
\end{proof}

For the Euclidean moment norm, the initial error is at most $\sqrt{|A_k|-1}\,e_k$, while the atom-projection error can be written directly as
$\Lip(\bv_{A_k})c_{\Lo}(C_He_t)^{\Lo/2}$. Thus the exponents are unchanged at fixed $k$ and fixed cliques. The maximum norm in \eqref{eq:maintransfer} avoids this coordinate-count conversion entirely. Although the kernel degree depends on the largest lifted clique, the constant and the global exponent $\Lo$ may depend on the whole constraint system. No ambient-dimension-uniform assertion follows without further geometry.

\section{Correlative sparsity with polynomial matrix constraints}\label{sec:matrix}
We extend the preceding argument to scalar objectives with clique-local polynomial matrix constraints. The matrix constraints change the local positivity tests but retain scalar sparse pseudo-moment sequences, so Lemmas~\ref{lem:gluing} and~\ref{lem:compactgluing} still apply. After recalling the matrix hierarchy, we bound the local negative-eigenvalue residual and use the same global projection.

\subsection{Scalar optimization with polynomial matrix constraints}\label{sec:matrixrecall}
We now consider the matrix description~\eqref{eq:matrixPOP} with prescribed clique-local blocks $G_{aj}\in\mathbb S^{p_{aj}}[\bx_{I_a}]$. Their inequalities are equivalent to $G(\bx)\succeq0$ for $G=\diag(G_{aj})$, but the sparse hierarchy keeps each block assigned to its own clique. The argument also allows additional clique-local equalities $h_{a\ell}=0$, imposed through their truncated ideals. After rescaling the bounding hypercube, write the feasible set as
\begin{equation}\label{eq:matrixset}
 \CX_G=\{\bx\in[0,1]^n:G_{aj}(\bx_{I_a})\succeq0,\
                         h_{a\ell}(\bx_{I_a})=0\ \text{for all }a,j,\ell\}.
\end{equation}
The description~\eqref{eq:matrixPOP} corresponds to an empty equality family. The objective $f=\sum_a f_a(\bx_{I_a})$ remains scalar. An entire matrix constraint is assigned to one variable clique: the union of the variable supports of all its entries must lie in that clique. This is the correlative pattern used for scalar polynomial optimization with matrix constraints by Nie et al.~\cite{niequ}.

In the matrix quadratic module of Scherer and Hol~\cite{schererhol}, a scalar certificate contains sums of $\bv^T G_{aj}\bv$ with polynomial vector multipliers. With the box generators included, its local truncation is
\begin{equation}\label{eq:matrixmodule}
 \QQ^G_{a,2r}
 =\left\{\sigma_0+\sum_{i\in I_a}\sigma_i b_i
          +\sum_{j,\ell}\bv_{j\ell}^TG_{aj}\bv_{j\ell}:
          \sigma_i\in\Sigma[\bx_{I_a}],\
          \text{all certificate degrees}\le2r\right\}.
\end{equation}
The degree of a vector-multiplier term is budgeted by
$2\deg \bv_{j\ell}+\deg G_{aj}$, without using cancellations. The dual uses scalar moments and the matrix localizer
\begin{equation}\label{eq:matrixlocalizer}
 M_{r-d_G}(G\biy)
 =\left(\sum_\delta G_\delta y_{\delta+\beta+\gamma}\right)_{
                   \beta,\gamma\in B_{r-d_G}^{I_a}},
 \qquad d_G=\lceil\deg G/2\rceil,
\end{equation}
where each displayed entry is a $p\times p$ matrix. Equivalently, its scalar rows are indexed by $(u,\beta)$, with $u\in[p]$ and $\beta\in B_{r-d_G}^{I_a}$. Positivity is equivalent to $\ell_{\biy}(\bv^TG\bv)\ge0$ for every admissible vector polynomial $\bv$. The moment matrices, equality conditions, normalization, and separator consistency are as before. We call each feasible array $\biy=(y_\alpha)_{\alpha\in A_{2r}}$ a sparse pseudo-moment sequence for the matrix-constrained hierarchy, with clique pseudo-moment restrictions $\biy^{(a)}$. We denote their feasible set by $\FF_r^{G,\QQ}$, its scalar objective infimum by $p_r^{G,\QQ}(f;\CX_G)$, and its truncated-moment distance by $\bd_{k,r}^{G,\QQ}(\CX_G)$.

Under running intersection, nonemptiness, and local Archimedean matrix quadratic modules, the scalar-objective sparse hierarchy converges; see Nie et al.~\cite[Section~1]{niequ}. This recalled convergence result supplies no quantitative exponent. We prove below a rate for the explicitly box-augmented hierarchy. Throughout this paper, every localizer in~\eqref{eq:matrixlocalizer} is imposed in full, including its off-diagonal matrix entries. The correlative graph concerns the variables in each entire matrix constraint. A matrix-valued objective leads to a different moment formulation; see Remark~\ref{rem:matrixscope}.

\subsection{A local kernel estimate for matrix feasibility}
For a symmetric matrix $A$, write $A_-$ for its negative part, with the negative eigenvalues replaced by their absolute values and the other eigenvalues by zero. Thus
\begin{displaymath}
 \|A_-\|_{\mathrm F}=\dist_{\mathrm F}(A,\mathbb S_+^p).
\end{displaymath}
For the clique-local description~\eqref{eq:matrixset}, use the block form of~\eqref{eq:matrixresidualrecall} and include the squared residuals of any additional equalities:
\begin{equation}\label{eq:matrixresidual}
 \rho_G(\bx)^2=\sum_{a,j}\|G_{aj}(\bx_{I_a})_-\|_{\mathrm F}^2
                  +\sum_{a,\ell}h_{a\ell}(\bx_{I_a})^2.
\end{equation}
For nonempty $\CX_G$, the matrix \L{}ojasiewicz inequality of Dinh and Pham~\cite[Theorem~4.1]{dinhpham}, recalled by Tran and Toh~\cite[Theorem~2.4]{tranmatrix}, provides the required distance bound. Apply it to $-\widetilde G$, where $\widetilde G$ is the block-diagonal polynomial matrix containing all $G_{aj}$, the scalar blocks $h_{a\ell}$ and $-h_{a\ell}$, and the box generators $b_i=x_i(1-x_i)$. Then $\CX_G=\{\bx:\widetilde G(\bx)\succeq0\}$. On $[0,1]^n$, its eigenvalue residual satisfies
\begin{displaymath}
\begin{aligned}
 \max\{0,-\lambda_{\min}(\widetilde G(\bx))\}
 &=\max\bigl(\{0\}\cup\{-\lambda_{\min}(G_{aj}(\bx_{I_a})):a,j\}
                  \cup\{|h_{a\ell}(\bx_{I_a})|:a,\ell\}\bigr)\\
 &\le\rho_G(\bx).
\end{aligned}
\end{displaymath}
Consequently, there exist $c_{\Lo}>0$ and $0<\Lo\le1$ such that
\begin{equation}\label{eq:matrixerrorbound}
 \dist(\bx,\CX_G)\le c_{\Lo}\rho_G(\bx)^{\Lo}
       \qquad(\bx\in[0,1]^n).
\end{equation}
Fix such a constant and exponent. This applies the cited matrix inequality to the whole constraint system; the exponent need not coincide with one obtained from a principal-minor description or from the individual cliques.

The following lemma applies to either of our kernels. Write
$\mathcal Kp(\bx)=\int p(\bt)K(\bt,\bx)d\varpi(\bt)$ and $M=\mathcal K1$ on one clique, where $\varpi$ denotes its product Chebyshev measure. The norm $\cheb{\cdot}$ is the shifted-Chebyshev coefficient $\ell^1$ norm defined in Section~\ref{sec:momentsets}.

\begin{lemma}\label{lem:matrixvariance}
Fix a clique $I_a$ and let $\biy^{(a)}=(y_\alpha^{(a)})_{\alpha\in B_{2r}^{I_a}}$ be a clique pseudo-moment sequence of order $r$, with $y_0^{(a)}=1$. Write $\ell=\ell_{\biy^{(a)}}$ for its Riesz functional, defined by
\begin{displaymath}
 \ell(p)=\ell_{\biy^{(a)}}(p)
 =\sum_{\alpha\in B_{2r}^{I_a}}p_\alpha y_\alpha^{(a)},
 \qquad p=\sum_{\alpha\in B_{2r}^{I_a}}p_\alpha\bx_{I_a}^{\alpha}
 \in\R[\bx_{I_a}]_{2r}.
\end{displaymath}
Let $G=(g_{uv})\in\mathbb S^p[\bx_{I_a}]$ have degree at most $d$. Suppose every polynomial evaluated by $\ell$ below has degree at most $2r$ and, for each $\bt$, the functional
$L_{\bt}(p)=\ell(K(\bt,\cdot)p)$ is nonnegative on squares of polynomials of degree at most $d$ and satisfies $L_{\bt}(G)\succeq0$. Suppose also that
\begin{displaymath}
\begin{gathered}
 \cheb{\mathcal Kp-p}\le b\cheb p\quad(\deg p\le2d),\\[3pt]
 |\ell(H)|\le\cheb H\quad\text{for }\quad
 H=\mathcal K(g_{uv}^2)-2g_{uv}\mathcal Kg_{uv}+g_{uv}^2M.
\end{gathered}
\end{displaymath}
Let $P(\bt)=L_{\bt}(1)$ and $Z=\ell(M)>0$. Then $d\mu=(P/Z)d\varpi$ is a probability measure and
\begin{equation}\label{eq:matrixvariance}
 \int\|G(\bt)_-\|_{\mathrm F}^2\,d\mu(\bt)
 \le \frac{4b}{Z}\sum_{u,v=1}^p\cheb{g_{uv}}^2.
\end{equation}
\end{lemma}
\begin{proof}
Positivity on squares gives $P\ge0$. For $P(\bt)>0$, the matrix
$B(\bt)=L_{\bt}(G)/P(\bt)$ is positive semidefinite. Cauchy--Schwarz for $L_{\bt}$ gives, entry by entry,
\begin{displaymath}
 P(\bt)|g_{uv}(\bt)-B_{uv}(\bt)|^2
 \le L_{\bt}\bigl((g_{uv}(\bt)-g_{uv}(\bx))^2\bigr).
\end{displaymath}
The distance of $G(\bt)$ to the PSD cone is no larger than its distance to $B(\bt)$. Summing these inequalities and integrating yields
\begin{displaymath}
 Z\int\|G(\bt)_-\|_{\mathrm F}^2\,d\mu(\bt)
 \le\sum_{u,v}\ell\left(\mathcal K(g_{uv}^2)
              -2g_{uv}\mathcal Kg_{uv}+g_{uv}^2M\right).
\end{displaymath}
The set $\{P=0\}$ has $\mu$-measure zero; the integrand on the right is nonnegative there as well. For an entry $g$, the polynomial in parentheses equals
\begin{displaymath}
 (\mathcal K(g^2)-g^2)-2g(\mathcal Kg-g)+g^2(M-1).
\end{displaymath}
Submultiplicativity of the coefficient norm bounds its norm by
$4b\cheb g^2$. The evaluation assumption proves~\eqref{eq:matrixvariance}.
\end{proof}

For the squared kernel, $K$ is SOS. Its product with a squared polynomial is therefore SOS, and $\ell(KG)\succeq0$ follows from the matrix localizer: for a constant vector $\bv$, expand $K=\sum_l q_l^2$ and evaluate $\sum_l(q_l\bv)^TG(q_l\bv)$. If
\begin{equation}\label{eq:matrixbudget}
 wD+2d\le r,
\end{equation}
all required positivity and evaluation bounds are valid. Indeed, the variance error polynomial has degree at most $wD+2d$, so Lemma~\ref{lem:evaluation} applies. The squared-kernel estimates~\eqref{eq:masscertificates}--\eqref{eq:squaredtensor} give
\begin{displaymath}
 b\le(1+\varepsilon_{2d})^w-1,\qquad Z\ge(1-\delta)^w.
\end{displaymath}
The fixed addition $2d$ in~\eqref{eq:matrixbudget} changes no asymptotic exponent. In particular, Lemma~\ref{lem:matrixvariance} gives a local squared matrix residual of order $(\log_2r)^3/r^2$.

\subsection{Correlative sparse rates with matrix constraints}
\begin{theorem}[Rate for the correlative sparse matrix quadratic module]\label{thm:matrixQ}
Let $\CX_G\ne\varnothing$, let the variable cliques satisfy running intersection, and impose the full local matrix localizers~\eqref{eq:matrixlocalizer}, explicit local box generators, and full separator consistency. Let $c_{\Lo}>0$ and $0<\Lo\le1$ be chosen as in the matrix \L{}ojasiewicz inequality~\eqref{eq:matrixerrorbound}. For fixed data and target degree $k$,
\begin{equation}\label{eq:matrixQrate}
 \bd_{k,r}^{G,\QQ}(\CX_G),\quad
 f_{\min,\CX_G}-p_r^{G,\QQ}(f;\CX_G)
 =O\left((\log_2r)^{3\Lo/2}/r^{\Lo}\right).
\end{equation}
The moment minimum equals the matching sparse SOS supremum. Both converge increasingly to $f_{\min,\CX_G}$.
\end{theorem}
\begin{proof}
Put $d=\max(\{1\}\cup\{\deg G_{aj}\}\cup\{\deg h_{a\ell}\})$. Choose a common squared kernel satisfying~\eqref{eq:matrixbudget}, with $q\ge\max\{k,2d\}$ and the logarithmic choice of $N$ in Corollary~\ref{thm:boxQ}. Let $\psi(r)=(\log_2r)^3/r^2$. The local probability measures have degree-$k$ moment error $O(\psi(r))$ and separator total-variation discrepancies at most $2w\delta=O(\psi(r))$.

By Lemma~\ref{lem:matrixvariance}, each local matrix residual has integral $O(\psi(r))$. Also $\ell_{\biy^{(a)}}(h_{a\ell}^2)=0$, so the local degree-$2d$ approximation gives $\int h_{a\ell}^2d\mu_a=O(\psi(r))$. Set
\begin{displaymath}
 \Phi_a=\sum_j\|G_{aj,-}\|_{\mathrm F}^2+\sum_\ell h_{a\ell}^2.
\end{displaymath}
These functions are continuous and nonnegative on their hypercubes, with fixed oscillations. Lemma~\ref{lem:compactgluing} therefore applies with $\eta=O(\psi(r))$. It gives the moment error
$O(\psi(r)+\psi(r)^{\Lo/2})$, which proves~\eqref{eq:matrixQrate}. Integrating the scalar objective gives the same bound for its optimal-value gap. This proof controls the residual of the glued measure before global projection; it does not require a polynomial square root of $G$ or a polynomial matrix slack substitution.

The box certificates $1\pm \bx^\alpha$ belong to the local matrix module. The compactness and interior-point separation argument of Proposition~\ref{prop:duality} thus applies without change. Restriction to a smaller order is feasible, proving monotonicity.
\end{proof}

With one clique the proof is the dense kernel-and-projection argument. Thus~\eqref{eq:matrixQrate} transfers its exponent to the correlative sparse hierarchy. Matrix sizes, clique widths, and coefficient norms enter the constants. The theorem does not identify $\Lo$ with an exponent belonging to a different dense penalty or residual construction.

\paragraph{A log-free variant with box-product localizers.}
The normalized kernel is positive in the box preordering, but need not be SOS. For matrix constraints, its use requires more than the usual localizers $M(G_{aj}\biy)$. Define the following explicit strengthening:
\begin{equation}\label{eq:matrixboxpreordering}
\begin{aligned}
 \mathcal B^G_{a,2r}
 &=\left\{\tau+\sum_{j,A,\ell}b_A\bv_{jA\ell}^TG_{aj}\bv_{jA\ell}:\right.\\
 &\hspace{12mm}\left.\tau\in\TT_{a,2r}([0,1]^{I_a}),\
       A\subseteq I_a,\ \text{all certificate degrees}\le2r\right\},\\
 b_A&=\prod_{i\in A}b_i.
\end{aligned}
\end{equation}
Equivalently, add matrix localizers for every $b_AG_{aj}$. This definition does not multiply two matrix generators and therefore does not require a convention for products of noncommuting matrices.

\begin{corollary}[Box-preordering-weighted matrix constraints]\label{cor:matrixB}
For the hierarchy defined by~\eqref{eq:matrixboxpreordering}, with the same equalities, consistency and error bound, the truncated-moment and scalar objective errors are $O(r^{-\Lo})$. The matching SOS supremum has the same rate.
\end{corollary}
\begin{proof}
Use $J_I=\prod_{i\in I}J_q$ with even $q$ and the sufficient budget $wq+2d\le r$. Its box-preordering certificate, multiplied by $\bv^TG\bv$, is tested by the added localizers. The functional $\ell(J_I\,\cdot)$ is also nonnegative on the required squares. Exact normalization gives $Z=1$ and equal separator marginals. The Jackson estimate~\eqref{eq:Jacksontensor} gives
$\cheb{\mathcal Kp-p}\le\Delta_{2d}(q)\cheb p$ for $\deg p\le2d$. Lemma~\ref{lem:matrixvariance} then bounds the local squared residual by $O(r^{-2})$. Apply compact gluing and projection as in Theorem~\ref{thm:matrixQ}. The same box interior argument gives moment--SOS equality.
\end{proof}

The log-free conclusion concerns~\eqref{eq:matrixboxpreordering}. It is not asserted for the unstrengthened correlative sparse matrix quadratic module. When the matrix generators are scalar, this cone is contained in the corresponding correlative sparse scalar preordering, so it is consistent with the scalar rate in Table~\ref{tab:rates}.

\begin{remark}[Scalar objectives and full matrix localizers]\label{rem:matrixscope}
Theorems in this section concern scalar sparse pseudo-moment sequences and their approximation by exact sparse moment sequences. Each polynomial matrix constraint is assigned in its entirety to one variable clique, and its full matrix localizer is imposed. Chordality of the variable graph does not justify replacing this localizer by selected principal blocks. Moreover, these results do not assert a rate for the original matrix-valued moment formulation of a matrix-valued objective. Miller et al.~\cite[Example~5.1]{miller} show that its correlative sparse hierarchy can fail to converge even under running intersection and local Archimedean assumptions. The scalar probability measures used here are essential to the gluing argument.
\end{remark}

\section{Conclusion}\label{sec:conclusion}
We have established convergence rates for moment--SOS hierarchies with correlative sparsity by constructing a joint approximating measure from local kernel measures. The quantitative gluing lemma controls separator discrepancies in total variation and preserves the dense hypercube rates for both the correlative sparse quadratic module and the correlative sparse preordering. The estimates apply uniformly to every sparse pseudo-moment sequence, and therefore control Hausdorff distances between truncated sparse pseudo-moment sets and exact sparse moment sets as well as polynomial objective gaps.

On compact sets, slack-variable lifting preserves the variable clique tree. Compact gluing controls the feasibility residual of the joint measure, and global atom projection converts this residual into the moment error through the exponent $\Lo$. For scalar objectives with polynomial matrix constraints, a kernel variance estimate supplies the corresponding residual bound directly. The global error-bound exponent and the constants depend on the constraint system; the gluing step itself causes no additional loss in the convergence exponent.

An exactly normalized kernel with quadratic-module certificates could remove the logarithmic factor while retaining separator compatibility. A further question is how quantitative regularity of the clique intersections controls the global feasibility exponent in terms of local geometry. Both questions concern the correlative sparse hierarchy studied here.

\appendix
\section{Auxiliary evaluation and duality facts}\label{app:auxiliary}
This appendix collects the elementary bounds and finite-order duality arguments used in the main text. All statements below use the explicit box generators specified in Section~\ref{sec:setup}.

\subsection{Evaluation bounds and correlative sparse duality}
\begin{lemma}[Evaluation and boundedness]\label{lem:evaluation}
Put $\tau_j(x)=T_j(2x-1)$ and $\tau_\alpha(\bx)=\prod_i\tau_{\alpha_i}(x_i)$. The coefficient $\ell^1$ norm is $\cheb p=\sum_\alpha|p_\alpha^{\mathrm{Cheb}}|$ for $p=\sum_\alpha p_\alpha^{\mathrm{Cheb}}\tau_\alpha$, as defined in Section~\ref{sec:momentsets}. If $\ell(1)=1$ and $\ell\ge0$ on the box quadratic module through degree $2r$, then
\begin{displaymath}
 |\ell(\bx^\alpha)|\le1\quad(|\alpha|\le2r),\qquad
 |\ell(p)|\le\cheb p\quad(\deg p\le r).
\end{displaymath}
The norm $\cheb{\cdot}$ is submultiplicative, and $\bx^\alpha$ has nonnegative shifted-Chebyshev coefficients summing to $1$.
\end{lemma}
\begin{proof}
The identities $1-x_i^2=(1-x_i)^2+2b_i$ and
\begin{displaymath}
 1-x_i^{2v}=(1-x_i^2)\sum_{j=0}^{v-1}x_i^{2j}
\end{displaymath}
give $1-\bx^{2\beta}\in\QQ_{2|\beta|}$. Split $\alpha=\beta+\gamma$ with $|\beta|,|\gamma|\le r$ and use
\begin{displaymath}
 1\pm \bx^\alpha=\tfrac12\bigl((\bx^\beta\pm \bx^\gamma)^2+1-\bx^{2\beta}+1-\bx^{2\gamma}\bigr).
\end{displaymath}
For Chebyshev evaluation, $1-\tau_j(x_i)^2=4b_iU_{j-1}(2x_i-1)^2$; telescoping gives $1-\tau_\alpha^2\in\QQ_{2|\alpha|}$. Positivity and Cauchy--Schwarz yield $|\ell(\tau_\alpha)|\le1$ for $|\alpha|\le r$. Expand $p$. Finally, $2\tau_i\tau_j=\tau_{i+j}+\tau_{|i-j|}$ proves submultiplicativity; $x=(1+\tau_1)/2$ proves the last assertion by repeated multiplication.
\end{proof}

\begin{proposition}[Moment and SOS values]\label{prop:duality}
The feasible set $\FF_r^{\mathcal C}(\CX)$ is compact, and strong duality holds, i.e.,
\begin{displaymath}
 p_r^{\mathcal C}(f;\CX)= s_r^{\mathcal C}(f;\CX).
 \end{displaymath}
\end{proposition}
\begin{proof}
Lemma~\ref{lem:evaluation} bounds every local coordinate, and the defining constraints are closed. Dirac measures on $\CX$ give feasibility. Let $E$ be the space of clique-supported polynomials of degree at most $2r$, and let $C=\mathcal C_{2r}^{\rm sp}(\CX)+\sum_a\JJ_{a,2r}$. The certificates $1\pm \bx^\alpha\in C$ show that $1$ is an interior point of $C$ in $E$: if $\sum_\alpha|a_\alpha|<1$, then $1+\sum_\alpha a_\alpha \bx^\alpha$ is a positive combination of these certificates and $1$.

Every nonzero functional nonnegative on $C$ has positive value at $1$, and can therefore be normalized to a sparse pseudo-moment functional. Separation consequently implies $f-p_r^{\mathcal C}\in\overline C$. For $\varepsilon>0$, adding the interior point $\varepsilon1$ gives $f-p_r^{\mathcal C}+\varepsilon\in C$, since $\operatorname{int}\overline C=\operatorname{int}C$ for a full-dimensional convex set. Weak duality gives the reverse inequality, proving equality of the two values.
\end{proof}

\section{Consistent local sequences and dense moment completion}\label{sec:gluingdistinction}
The distinction between measure gluing and completion of sparse moments is important for the rate-transfer argument. We recall the relevant results from the literature, separating consistency of a sparse functional, completion of a positive semidefinite matrix, and extension to a representing measure.

\paragraph{Consistent sparse functionals.}
In the correlative sparse hierarchy of Lasserre~\cite[Section~3]{lasserre2006}, shared monomials have a single moment coordinate. Thus the consistency conditions~\eqref{eq:consistency} identify the clique restrictions with one functional on $E_{2r}=\sum_a\R[\bx_{I_a}]_{2r}$, given by
\begin{equation}\label{eq:sparsefunctional}
 \ell\left(\sum_a p_a\right)=\sum_a\ell_{\biy^{(a)}}(p_a),
 \qquad p_a\in\R[\bx_{I_a}]_{2r}.
\end{equation}
This is the sparse pseudo-moment functional used in the paper. Its domain contains only clique-supported polynomials; a dense extension must also assign the missing mixed moments and satisfy the dense positivity conditions.

\paragraph{Matrix completion and moment extension.}
The chordal completion theorem of Grone et al.~\cite{grone} characterizes positive semidefinite completability of a partial matrix by positivity of its clique principal submatrices; see also Vandenberghe and Andersen~\cite[Theorem~10.1]{vandenberghe}. A moment matrix has the additional Hankel identities $M_{\alpha,\beta}=M_{\gamma,\delta}$ whenever $\alpha+\beta=\gamma+\delta$, and all localizing matrices must use the same moment coordinates. These requirements are not part of the chordal matrix-completion theorem. In particular, chordality of the variable graph alone does not supply a dense moment completion.

Moment-extension results impose additional hypotheses. Laurent and Mourrain~\cite[Theorem~1.4]{laurentmourrain} prove a flat extension theorem for a finite monomial set $C$ connected to $1$, with $C^+=C\cup\bigcup_i x_iC$: the condition $\operatorname{rank}M_{C^+}=\operatorname{rank}M_C$ gives a unique infinite flat moment extension. Under positivity, their Theorem~3.2 gives a representing measure. These flatness and positivity hypotheses are stronger than agreement of truncated separator moments and are not assumed in our convergence estimates.

\paragraph{Gluing probability measures.}
Lasserre~\cite[Section~6.1 and Lemma~6.4]{lasserre2006} proves that, under running intersection, local probability measures with equal separator marginals admit a joint probability measure, with the corresponding support constraints preserved. The hypothesis concerns equality of measures. Equality of all separator moments suffices on compact sets, but equality through a fixed degree need not determine a marginal measure. This distinction is also reflected in Tchakaloff's theorem: a probability measure on a compact set has a finitely atomic representative of any fixed finite collection of polynomial moments; see Bayer and Teichmann~\cite[Corollary~2]{bayer}.

Our construction uses these measure-level requirements rather than assuming an exact dense completion of the original sparse pseudo-moment sequence. Lemma~\ref{lem:compatiblekernel} produces local measures with identical separator marginals, while Lemma~\ref{lem:squaredseparator} controls their total-variation discrepancies. The joint measure then approximates the prescribed truncated moments. This approximation does not assert preservation of every original sparse coordinate by a dense pseudo-moment extension at the same order.

\begingroup
\small
\nocite{*}
\bibliographystyle{plain}
\bibliography{sections/references}

@article{baldi,
  author  = {Baldi, L. and Mourrain, B.},
  title   = {{On the effective Putinar's Positivstellensatz and moment approximation}},
  journal = {Mathematical Programming},
  volume  = {200},
  year    = {2023},
  pages   = {71--103}
}

@article{bayer,
  author  = {Bayer, C. and Teichmann, J.},
  title   = {{The proof of Tchakaloff's theorem}},
  journal = {Proceedings of the American Mathematical Society},
  volume  = {134},
  year    = {2006},
  pages   = {3035--3040}
}

@misc{optimal2026,
  author       = {Blomenhofer, A. T. and Gribling, S. and Henrion, D. and Korda, M. and Kumar, P. and Magron, V. and Slot, L. and Vargas, L. F. and Vera, J.},
  title        = {{Optimal degree bounds for sum-of-squares representations}},
  howpublished = {arXiv:2610.24638v1},
  year         = {2026},
  eprint       = {2610.24638},
  archivePrefix = {arXiv},
  url          = {https://arxiv.org/abs/2610.24638},
  note         = {\url{https://arxiv.org/abs/2610.24638}}
}

@book{bcr,
  author    = {Bochnak, J. and Coste, M. and Roy, M.-F.},
  title     = {{Real Algebraic Geometry}},
  series    = {Ergebnisse der Mathematik und ihrer Grenzgebiete},
  volume    = {36},
  publisher = {Springer},
  year      = {1998}
}

@article{changpollard,
  author  = {Chang, J. T. and Pollard, D.},
  title   = {{Conditioning as disintegration}},
  journal = {Statistica Neerlandica},
  volume  = {51},
  year    = {1997},
  pages   = {287--317},
  doi     = {10.1111/1467-9574.00056},
  url     = {https://doi.org/10.1111/1467-9574.00056},
  note    = {\url{https://doi.org/10.1111/1467-9574.00056}}
}

@misc{kernel,
  author       = {Gribling, S. and de Klerk, E. and Vera, J. C.},
  title        = {{Squared polynomial approximation kernels for the hypercube: improved error bounds and implications for Lasserre hierarchies}},
  howpublished = {arXiv:2605.31496v1},
  year         = {2026},
  eprint       = {2605.31496},
  archivePrefix = {arXiv},
  url          = {https://arxiv.org/abs/2605.31496},
  note         = {\url{https://arxiv.org/abs/2605.31496}}
}

@article{grone,
  author  = {Grone, R. and Johnson, C. R. and S\'a, E. M. and Wolkowicz, H.},
  title   = {{Positive definite completions of partial Hermitian matrices}},
  journal = {Linear Algebra and its Applications},
  volume  = {58},
  year    = {1984},
  pages   = {109--124},
  doi     = {10.1016/0024-3795(84)90207-6},
  url     = {https://doi.org/10.1016/0024-3795(84)90207-6},
  note    = {\url{https://doi.org/10.1016/0024-3795(84)90207-6}}
}

@article{grimm,
  author  = {Grimm, D. and Netzer, T. and Schweighofer, M.},
  title   = {{A note on the representation of positive polynomials with structured sparsity}},
  journal = {Archiv der Mathematik},
  volume  = {89},
  year    = {2007},
  pages   = {399--403},
  doi     = {10.1007/s00013-007-2234-z},
  url     = {https://doi.org/10.1007/s00013-007-2234-z},
  note    = {\url{https://doi.org/10.1007/s00013-007-2234-z}}
}

@book{kallenberg,
  author    = {Kallenberg, O.},
  title     = {{Foundations of Modern Probability}},
  edition   = {Third},
  series    = {Probability Theory and Stochastic Modelling},
  volume    = {99},
  publisher = {Springer},
  year      = {2021},
  doi       = {10.1007/978-3-030-61871-1},
  url       = {https://doi.org/10.1007/978-3-030-61871-1},
  note      = {\url{https://doi.org/10.1007/978-3-030-61871-1}}
}

@article{korda,
  author  = {Korda, M. and Magron, V. and R\'ios-Zertuche, R.},
  title   = {{Convergence rates for sums-of-squares hierarchies with correlative sparsity}},
  journal = {Mathematical Programming},
  volume  = {209},
  year    = {2025},
  pages   = {435--473},
  doi     = {10.1007/s10107-024-02071-6},
  url     = {https://doi.org/10.1007/s10107-024-02071-6},
  note    = {\url{https://doi.org/10.1007/s10107-024-02071-6}}
}

@article{lasserre2001,
  author  = {Lasserre, J.-B.},
  title   = {{Global optimization with polynomials and the problem of moments}},
  journal = {SIAM Journal on Optimization},
  volume  = {11},
  year    = {2001},
  pages   = {796--817}
}

@article{lasserre2006,
  author  = {Lasserre, J.-B.},
  title   = {{Convergent SDP-relaxations in polynomial optimization with sparsity}},
  journal = {SIAM Journal on Optimization},
  volume  = {17},
  year    = {2006},
  pages   = {822--843}
}

@article{laurent,
  author  = {Laurent, M. and Slot, L.},
  title   = {{An effective version of Schm\"udgen's Positivstellensatz for the hypercube}},
  journal = {Optimization Letters},
  volume  = {17},
  year    = {2023},
  pages   = {515--530},
  doi     = {10.1007/s11590-022-01922-5},
  url     = {https://doi.org/10.1007/s11590-022-01922-5},
  note    = {\url{https://doi.org/10.1007/s11590-022-01922-5}}
}

@book{marshall,
  author    = {Marshall, M.},
  title     = {{Positive Polynomials and Sums of Squares}},
  series    = {Mathematical Surveys and Monographs},
  volume    = {146},
  publisher = {American Mathematical Society},
  year      = {2008}
}

@phdthesis{parrilo,
  author = {Parrilo, P. A.},
  title  = {{Structured semidefinite programs and semialgebraic geometry methods in robustness and optimization}},
  school = {California Institute of Technology},
  year   = {2000}
}

@article{putinar,
  author  = {Putinar, M.},
  title   = {{Positive polynomials on compact semi-algebraic sets}},
  journal = {Indiana University Mathematics Journal},
  volume  = {42},
  year    = {1993},
  pages   = {969--984}
}

@article{schmudgen,
  author  = {Schm\"udgen, K.},
  title   = {{The K-moment problem for compact semi-algebraic sets}},
  journal = {Mathematische Annalen},
  volume  = {289},
  year    = {1991},
  pages   = {203--206}
}

@article{tran,
  author  = {Tran, H. A. and Toh, K.-C.},
  title   = {{On the convergence rates of moment-SOS hierarchies approximation of truncated moment sequences}},
  journal = {Mathematical Programming},
  year    = {2026},
  doi     = {10.1007/s10107-026-02394-6},
  url     = {https://doi.org/10.1007/s10107-026-02394-6},
  note    = {Published online 7 July 2026. \url{https://doi.org/10.1007/s10107-026-02394-6}}
}

@unpublished{extension,
  author = {Tran, H. A. and Toh, K.-C.},
  title  = {{Extension on the convergence rates of moment-SOS hierarchies via approximation of truncated moment sequences}},
  year   = {2026},
  note   = {Unpublished manuscript}
}

@article{tranmatrix,
  author  = {Tran, H. A. and Toh, K.-C.},
  title   = {{Convergence rates of sum-of-squares hierarchies for polynomial semidefinite programs}},
  journal = {SIAM Journal on Optimization},
  volume  = {36},
  year    = {2026},
  pages   = {760--790},
  doi     = {10.1137/24M1670184},
  url     = {https://doi.org/10.1137/24M1670184},
  note    = {\url{https://doi.org/10.1137/24M1670184}}
}

@article{waki,
  author  = {Waki, H. and Kim, S. and Kojima, M. and Muramatsu, M.},
  title   = {{Sums of squares and semidefinite program relaxations for polynomial optimization problems with structured sparsity}},
  journal = {SIAM Journal on Optimization},
  volume  = {17},
  year    = {2006},
  pages   = {218--242}
}

@article{tssos,
  author  = {Wang, J. and Magron, V. and Lasserre, J.-B.},
  title   = {{TSSOS: A moment-SOS hierarchy that exploits term sparsity}},
  journal = {SIAM Journal on Optimization},
  volume  = {31},
  year    = {2021},
  pages   = {30--58},
  url     = {https://arxiv.org/abs/1912.08899},
  note    = {\url{https://arxiv.org/abs/1912.08899}}
}

@article{chordal,
  author  = {Wang, J. and Magron, V. and Lasserre, J.-B.},
  title   = {{Chordal-TSSOS: A moment-SOS hierarchy that exploits term sparsity with chordal extension}},
  journal = {SIAM Journal on Optimization},
  volume  = {31},
  year    = {2021},
  pages   = {114--141},
  doi     = {10.1137/20M1323564},
  url     = {https://doi.org/10.1137/20M1323564},
  note    = {\url{https://doi.org/10.1137/20M1323564}}
}

@article{cstssos,
  author  = {Wang, J. and Magron, V. and Lasserre, J.-B. and Mai, N. H. A.},
  title   = {{CS-TSSOS: Correlative and term sparsity for large-scale polynomial optimization}},
  journal = {ACM Transactions on Mathematical Software},
  volume  = {48},
  year    = {2022},
  pages   = {1--26},
  eid     = {42},
  doi     = {10.1145/3569709},
  url     = {https://doi.org/10.1145/3569709},
  note    = {Article 42. \url{https://doi.org/10.1145/3569709}}
}

@article{niequ,
  author  = {Nie, J. and Qu, Z. and Tang, X. and Zhang, L.},
  title   = {{Sparse polynomial optimization with matrix constraints}},
  journal = {Journal of Global Optimization},
  volume  = {95},
  year    = {2026},
  pages   = {711--739},
  doi     = {10.1007/s10898-025-01539-9},
  url     = {https://doi.org/10.1007/s10898-025-01539-9},
  note    = {\url{https://doi.org/10.1007/s10898-025-01539-9}}
}

@article{miller,
  author  = {Miller, J. and Wang, J. and Guo, F.},
  title   = {{Sparse polynomial matrix optimization}},
  journal = {SIAM Journal on Optimization},
  volume  = {36},
  year    = {2026},
  pages   = {503--533},
  doi     = {10.1137/24M1719761},
  url     = {https://doi.org/10.1137/24M1719761},
  note    = {\url{https://doi.org/10.1137/24M1719761}}
}

@article{schererhol,
  author  = {Scherer, C. W. and Hol, C. W. J.},
  title   = {{Matrix sum-of-squares relaxations for robust semi-definite programs}},
  journal = {Mathematical Programming},
  volume  = {107},
  year    = {2006},
  pages   = {189--211},
  doi     = {10.1007/s10107-005-0684-2},
  url     = {https://doi.org/10.1007/s10107-005-0684-2},
  note    = {\url{https://doi.org/10.1007/s10107-005-0684-2}}
}

@article{laurentmourrain,
  author  = {Laurent, M. and Mourrain, B.},
  title   = {{A generalized flat extension theorem for moment matrices}},
  journal = {Archiv der Mathematik},
  volume  = {93},
  number  = {1},
  pages   = {87--98},
  year    = {2009},
  doi     = {10.1007/s00013-009-0007-6}
}

@article{vandenberghe,
  author  = {Vandenberghe, L. and Andersen, M. S.},
  title   = {{Chordal graphs and semidefinite optimization}},
  journal = {Foundations and Trends in Optimization},
  volume  = {1},
  number  = {4},
  pages   = {241--433},
  year    = {2015},
  doi     = {10.1561/2400000006}
}

@article{bmp,
  author  = {Baldi, L. and Mourrain, B. and Parusi\'nski, A.},
  title   = {{On {\L}ojasiewicz inequalities and the effective Putinar's Positivstellensatz}},
  journal = {Journal of Algebra},
  volume  = {662},
  pages   = {741--767},
  year    = {2025},
  doi     = {10.1016/j.jalgebra.2024.08.022}
}

@article{dinhpham,
  author  = {Dinh, S. T. and Pham, T. S.},
  title   = {{{\L}ojasiewicz-type inequalities with explicit exponents for the largest eigenvalue function of real symmetric polynomial matrices}},
  journal = {International Journal of Mathematics},
  volume  = {27},
  number  = {2},
  pages   = {1650012},
  year    = {2016},
  doi     = {10.1142/S0129167X16500129}
}
\endgroup

\end{document}